\documentclass[10pt]{article}

\usepackage{color}
\usepackage{latexsym,dsfont}
\usepackage{amssymb}
\usepackage{amsmath}
\usepackage{amsthm}
\usepackage{graphicx}
\usepackage{enumerate}
\usepackage{hyperref}

\newtheorem{Theorem}{Theorem}[part]
\newtheorem{Definition}{Definition}[part]
\newtheorem{Proposition}{Proposition}[part]

\newtheorem{Lemma}{Lemma}[part]

\newtheorem{Corollary}{Corollary}[part]
\newtheorem{Remark}{Remark}[part]

\def \R{\mathbb{R}}

\def \E{\mathbb{E}}
\def \F{\mathbb{F}}
\def \G{\mathbb{G}}

\def \P{\mathbb{P}}
\def \Q{\mathbb{Q}}
\def \D{\mathbb{D}}

\def \PMc{{\cal P \cal M}}

\def \Bc{{\cal B}}

\def \Dc{{\cal D}}
\def \Ec{{\cal E}}
\def \Fc{{\cal F}}
\def \Gc{{\cal G}}

\def \Pc{{\cal P}}
\def \Mc{{\cal M}}

\def \Sc{{\cal S}}
\def \Tc{{\cal T}}

\def \1{\mathds{1}}
\def \XunPi{X^{1, \pi}}
\def \XtilPi{\tilde{X}^\pi}
\def \YunPi{Y^{1, \pi}}
\def \ZunPi{Z^{1, \pi}}

\def \YzePi{Y^{0, \pi}}
\def \ZzePi{Z^{0, \pi}}
\def \uZzePi{\underline Z^{0, \pi}}
\def \YtilzePi{\tilde Y^{0, \pi}}
\def \ZtilzePi{\tilde Z^{0, \pi}}
\def \uZtilzePi{\underline{\tilde Z}^{0, \pi}}
\def \XzePi{X^{0, \pi}}

\def \ni{\noindent}

\def \eps{\varepsilon}

\def \ep{\hbox{ }\hfill$\Box$}

\def\reff#1{{\rm(\ref{#1})}}

\def\beqs{\begin{eqnarray*}}
\def\enqs{\end{eqnarray*}}
\def\beq{\begin{eqnarray}}
\def\enq{\end{eqnarray}}

\newcommand{\nc}{\newcommand}
\nc{\esssup}{\mathop{\mathrm{ess\;sup}}}

\begin{document}
\title{A decomposition approach for\\ the discrete-time approximation of \\FBSDEs with a jump}
\author{Idris Kharroubi \footnote{The research of the author benefited from the support of the French ANR research grant LIQUIRISK} \\
\small{CEREMADE, CNRS  UMR 7534}\\
\small{Universit\'e Paris Dauphine} \\
\small{\texttt{kharroubi @ ceremade.dauphine.fr}}\and 
Thomas Lim \footnote{The research of the author benefited from the support of the ``Chaire March\'e en mutation'', F\'ed\'eration Bancaire Fran\c caise}\\\small{ENSIIE,}\\
\small{Laboratoire de Math\'ematiques et Mod\'elisation d'Evry,}\\
\small{  CNRS  UMR 8071} \\
\small{\texttt{thomas.lim @ ensiie.fr} }}
\date{}


\maketitle
\begin{abstract}
We are concerned with the discretization of a solution of a Forward-Backward stochastic differential equation (FBSDE) with a jump process depending on the Brownian motion. In this paper, we study the cases of Lipschitz generators and the generators with a quadratic growth w.r.t. the variable $z$. We propose a recursive scheme based on a general existence result given in the companion paper \cite{kl11} and we study the error induced by the time discretization. We prove the convergence of the scheme when the number of time steps $n$ goes to infinity. Our approach allows to get a convergence rate similar to that of schemes of Brownian FBSDEs. 

\end{abstract}

\vspace{1cm}

\ni\textbf{Keywords:} discrete-time approximation, forward-backward SDE, Lipschitz generator, generator of quadratic growth, 
progressive enlargement of filtrations, decomposition in the reference filtration. \\

\vspace{1cm}

\ni\textbf{MSC classification (2010):} 65C99, 60J75, 60G57.


\setcounter{section}{0}

\section{Introduction}
\setcounter{equation}{0} \setcounter{Assumption}{0}
\setcounter{Theorem}{0} \setcounter{Proposition}{0}
\setcounter{Corollary}{0} \setcounter{Lemma}{0}
\setcounter{Definition}{0} \setcounter{Remark}{0}

In this paper, we study a discrete-time approximation for the solution of a forward-backward stochastic differential equation (FBSDE) with a jump of the form
 \begin{equation*} \left\{
\begin{aligned}
  X_t ~ & = ~ x + \int_0^t b( s, X_s) ds + \int_0^t \sigma(s, X_s) dW_s + \int_0^t \beta(s, X_{s^-}) dH_s \;,\\
  Y_t ~ & = ~ g(X_T) + \int_t^T f( s, X_s, Y_s, Z_s, U_s) ds - \int_t^T Z_s dW_s - \int_t^T U_s dH_s\;, \\
  \end{aligned}
\right.\end{equation*}
where $H_t = \1_{\tau \leq t}$ and $\tau$ is a jump time, which can represent a default time in credit risk or counterparty risk. Such equations naturally appear in finance, see for example Bielecki and Jeanblanc \cite{biejea08}, Lim and Quenez \cite{limque09}, Peng and Xu \cite{penxu10}, Ankirchner \emph{et al.} \cite{ankblaeyr09} for an application to exponential utility maximization problem and Kharroubi and Lim \cite{kl11} for the hedging problem in a complete market. The approximation of such equation is therefore of important interest for practical applications in finance. In this paper, we study the case where the generator $f$ is Lipschitz or with quadratic growth w.r.t. $Z$. 

In the literature, the problem of discretization of FBSDEs  with Lipschitz generator has been widely studied in the Brownian framework, \textit{i.e.} no jump, see e.g. \cite{maproyon94, doumapro96, che97, boutou04, 
zha04, delmen06}. 
More recently, the case of quadratic generators w.r.t. $Z$ has been considered by Imkeller \emph{et al.} \cite{imkreizha10} and Richou \cite{ric10}.
For Lipschitz generators, the discrete-time approximation of FBSDEs with jumps is studied by Bouchard and Elie \cite{boueli08} in the case of Poissonian jumps independent of the Brownian motion. Their approach is based on a regularity result for the process $Z$, which is given by Malliavin calculus tools. This regularity result for the process $Z$ was first proved by Zhang \cite{zha04} in a Brownian framework to provide a convergence rate for the discrete-time approximation of FBSDEs. The use of Malliavin calculus to prove regularity on $Z$ is possible in \cite{boueli08} since the authors suppose that the Brownian motion is independent of the jump measure. 

In our case,  we only assume that the random jump time $\tau$ admits a conditional density given $W$, which is assumed to be absolutely continuous w.r.t. the Lebesgue measure. In particular, we do not specify a particular law for $\tau$ and we do not assume that $\tau$ is independent of $W$ as for the case of a Poisson random measure. 

To the best of our knowledge, no Malliavin calculus theory has been set for such a framework. 
Thus, the method used in \cite{boueli08} fails to provide a convergence rate for the approximation in this context. 

We therefore follow another approach, which consists in using the decomposition result given in the companion paper \cite{kl11} to write the solution of a FBSDE with a jump as a combination of solutions to a recursive system of FBSDEs without jump.
We then prove a regularity result on the $Z$ components of Brownian BSDEs coming from the decomposition of the BSDE with a jump.  This regularity result allows to get a rate for the convergence of the discrete-time schemes for these BSDEs as in \cite{zha04} or \cite{boueli08} for the Lipschitz case and \cite{ric10} for the quadratic case.

Finally, we recombine the approximations of the solutions to recursive system of Brownian FBSDEs to get a discretization of the solution to the FBSDE with a jump. 

We notice that our approach also allows to weaken the assumption on the forward jump coefficient in the Lipschitz case. More precisely, we only assume that $\beta$ is Lipschitz continuous, unlike \cite{boueli08}  supposing that $\beta$ is regular and  the matrix $I_d+\nabla\beta$ is elliptic.



 As said above, this kind of FBSDEs with a jump appears in finance. The general assumptions made on the jump time $\tau$ allow to modelize  general phenomena as a firm default or simpler as  a jump of an asset that can be seen as contagion from the default of another firm on the market, see e.g. \cite{jkp11} for some examples. In particular, the approximation of these FBSDEs has its own interest, since it provides approximations of optimal gains and strategies of the studied investment problems.
   
  We choose to present our results in the case of a single jump and a one-dimensional Brownian motion for the sake of simplicity. We notice that they can easily be extended to the case of a $d$-dimensional Brownian motion and multiple jumps with eventually random marks, as in \cite{kl11}, taking values in a finite space.
   
The paper is organized as follows. The next section presents the framework of progressive enlargement of a Brownian filtration by a random jump, and the well posedness of FBSDEs in this context.  In Section 3, we present the discrete-time schemes for the forward and backward solutions based on the decomposition given in the previous section. Finally, in Section 4, we study  the convergence  rate of the scheme for the forward solution.  In Sections 5 and 6, we study  the convergence  rate of the scheme for the backward solution for the Lipschitz case respectively for the quadratic case. 

\section{Preliminaries}
\setcounter{equation}{0} \setcounter{Assumption}{0}
\setcounter{Theorem}{0} \setcounter{Proposition}{0}
\setcounter{Corollary}{0} \setcounter{Lemma}{0}
\setcounter{Definition}{0} \setcounter{Remark}{0}

\subsection{Notation}
Throughout this paper, we let $(\Omega, \Gc, \P)$ a complete probability space on which is defined a standard one dimensional Brownian motion $W$. We denote 
 $\F :={(\Fc_t)}_{t\geq 0}$ the natural filtration of $W$ augmented by all the $\P$-null sets. We also consider on this space a random time $\tau$, i.e. a nonnegative $\Fc$-measurable random variable, and we denote classically the associated jump process by $H$ which is given by 
\beqs
H_t & := & \1_{\tau \leq t}\;,\quad t\geq 0\;. 
\enqs
We denote by  $\D := (\Dc_t)_{t \geq 0}$ the smallest right-continuous filtration for which $\tau$ is a stopping time. The global information is then defined by the progressive enlargement $\G := {(\Gc_t)}_{t \geq 0}$ of the initial filtration where \beqs
\Gc_t &  :=  & \bigcap_{\eps>0} \Big(\Fc_{t+\eps} \vee \Dc_{t+\eps}\Big)
\enqs
for all $t\geq 0$. 
 This kind of enlargement  was introduced by Jacod, Jeulin and Yor in the 80s (see e.g. \cite{jeu80}, \cite{jeuyor85} and \cite{jac87}). We introduce some notations used throughout the paper\begin{itemize}
\item $\Pc(\F)$ (resp. $\Pc(\G)$) is the $\sigma$-algebra of $\F$ (resp. $\G$)-predictable measurable subsets of $\Omega \times \R_{+}$, i.e. the $\sigma$-algebra generated by the left-continuous $\F$ (resp. $\G$)-adapted processes, 
\item $\PMc(\F)$ (resp. $\PMc(\G)$) is the $\sigma$-algebra of $\F$ (resp. $\G$)-progressively measurable subsets of $\Omega\times\R_{+}$. 
\end{itemize}
We shall make, throughout the sequel, the standing assumption in the progressive enlargement of filtrations known as  density assumption (see e.g. \cite{jiapha09, jkp11,kl11}).

\vspace{2mm}

\ni \textbf{(DH)} There exists a positive and bounded $\Pc(\F)\otimes\Bc(\R_+)$-measurable process $\gamma$ such that
\beqs
\P \big[\tau\in d\theta \;\big| \;\Fc_t \big] & = & \gamma_t(\theta)d\theta \;,\quad t\geq 0\;.
\enqs

\vspace{2mm}

Using Proposition 2.1 in \cite{kl11} we get that \textbf{(DH)}  ensures that the process $H$ admits an intensity.
 \begin{Proposition}
The process $H$ admits a compensator of the form $\lambda_tdt$, where the process $\lambda$ is defined by 
\beqs
\lambda_t & := & \frac{\gamma_t(t)}{\P \big[\tau>t\; \big| \;\Fc_t\big]}\mathds{1}_{t\leq \tau}\;,\quad  t\geq 0\;.
\enqs
\end{Proposition} 
\ni We impose the following assumption to the process $\lambda$.

\vspace{2mm}

\ni\textbf{(HBI)} The process $\lambda$ is bounded. 

\vspace{2mm}

\ni We also introduce the martingale invariance assumption known as the \textbf{(H)}-hypothesis.

\vspace{2mm}

\ni \textbf{(H)} Any $\F$-martingale remains a $\G$-martingale.

\vspace{2mm}

 
\ni We now introduce the following spaces, where $a,b\in\R_+$ with $a< b$,
and $T < \infty$ is the terminal time.
\begin{itemize}
\item $\Sc_{\G}^\infty[a,b]$ (resp. $\Sc_{\F}^\infty[a,b]$) is the set of $\PMc(\G)$ (resp.  $\PMc(\F)$)-measurable processes $(Y_t)_{t\in[a,b]}$ essentially bounded
\beqs
{\| Y \|}_{\Sc^\infty[a,b]} & := & \esssup_{t\in[a,b]}|Y_{t}|~<~\infty \;.
\enqs
\item $\Sc_{\G}^p[a,b]$ (resp. $\Sc_{\F}^p[a,b]$), with $p\geq 2$, is the set of $\PMc(\G)$ (resp.  $\PMc(\F)$)-measurable processes $(Y_t)_{t\in[a,b]}$ such that
\beqs
{\| Y \|}_{\Sc^p[a,b]} & := & \Big(\E\Big[\sup_{t\in[a,b]}|Y_{t}|^p\Big]\Big)^{1\over p}~<~\infty \;.
\enqs

\item $H^p_{\G}[a,b]$ (resp. $H^p_{\F}[a,b]$), with $p\geq 2$, is the set of $\Pc(\G)$ (resp. $\Pc(\F)$)-measurable processes $(Z_t)_{t\in[a,b]}$ such that
\beqs
{\|Z\|}_{H^p[a,b]}  & := & \E\Big[ \Big(\int_a^b |Z_t|^2 dt \Big)^{p\over 2}\Big]^{1\over p} ~< ~\infty \;.
\enqs
\item $L^2(\lambda)$ is the set of $\Pc(\G)$-measurable processes $(U_t)_{t \in [0, T]}$ such that
\beqs
{\|U\|}_{L^2(\mu)} &:= & \Big(\E\Big[\int_0^T \lambda_s |U_{s}|^2ds\Big]\Big)^{1\over2}~<~\infty\;.
\enqs

\end{itemize}


  \subsection{ Forward-Backward SDE with a jump}
Given measurable functions $b:[0,T]\times\R \rightarrow \R$, $\sigma:[0,T]\times\R \rightarrow \R$, $\beta: [0,T]\times\R \rightarrow\R$, $g: \R \rightarrow\R$ and $f: [0, T] \times \R \times\R\times\R \times \R\rightarrow\R$, and an initial condition $x \in \R$, we  study the discrete-time approximation of  the solution $(X, Y, Z, U)$ in $\Sc^2_\G [0, T] \times \Sc^\infty_\G [0, T] \times H^2_\G [0, T] \times L^2(\lambda)$ to the following forward-backward stochastic differential equation
\begin{eqnarray}\label{FSDE} 
  X_t  & = & x + \int_0^t \hspace{-2mm}b(s, X_s) ds + \int_0^t \hspace{-2mm}\sigma(s, X_s) dW_s + \int_0^t\hspace{-2mm}\beta(s, X_{s^-}) dH_s \;, \quad 0 \leq t\leq T\;,\\ \nonumber
  Y_t  & =  & g(X_T) + \int_t^T \hspace{-2mm}f\big(s, X_s, Y_s, Z_s, (1-H_s)U_s\big) ds\\ \label{BSDE}
   & & \quad \quad \quad \quad  - \int_t^T\hspace{-2mm} Z_s dW_s - \int_t^T\hspace{-2mm} U_s dH_s\;, \quad 0 \leq t\leq T\;, ~~
\end{eqnarray}
  when the generator of the BSDE is Lipschitz or has a quadratic growth w.r.t. $Z$. 
  \begin{Remark}
  {\rm In the BSDE \reff{BSDE}, the jump component $U$ of the unknown $(Y,Z,U)$ appears in the generator $f$ with the additional multiplicative term $1-H$. This ensures the equation to be well posed in $\Sc_\G^\infty[0,T]\times H^2_\G[0,T]\times L^2(\lambda)$. Indeed,  the component $U$ lives in $L^2(\lambda)$, thus its value on $(\tau\wedge T,T]$ is not defined since the intensity $\lambda$ vanishes on $(\tau\wedge T,T]$. We therefore introduce the term $1-H$ to kill the value of $U$ on $(\tau\wedge T,T]$ and hence to avoid making the equation depending on it.   }
  \end{Remark}
  We first prove that the decoupled  system \reff{FSDE}-\reff{BSDE} admits a solution.
  To this end, we introduce several  assumptions on the coefficients  $b$, $\sigma$, $\beta$, $g$ and $f$. 
  We consider the following assumption for the forward coefficients.  
  
\vspace{2mm}  
  
 \ni\textbf{(HF)} There exists a constant $K$ such that 
  the functions $b$, $\sigma$ and $\beta$ satisfy 
  \beqs
  |b(t, 0)|+ |\sigma (t,0 )|+|\beta(t, 0)| & \leq  & K \;,
  \enqs
and
  \beqs
  |b(t, x) - b(t, x')| + |\sigma (t, x) - \sigma (t, x')| + |\beta(t, x) - \beta(t, x')| & \leq & K |x - x'| \;,
  \enqs
  for all $(t,x,x')\in[0,T]\times\R \times\R$. 
  
\vspace{2mm}  
  

  

\ni For the backward coefficients $g$ and $f$, we impose the following assumptions for the Lipschitz case.
\vspace{2mm}

 \ni\textbf{(HBL)} There exists a constant $K$ such that the functions $g$ and $f$ satisfy
  \beqs
  |f(t,x,0,0,0)| + |g(x)| & \leq & K\;,
  \enqs
  and
  \beqs
  |f (t, x,y,z,u) - f ( t,x,y',z',u')|  & \leq & 
  K \big( |y-y'|+|z-z'|+|u-u'|\big) \;,
  \enqs
  for all $(t,x,y,y',z,z',u,u')\in[0,T]\times\R \times\R^2\times\R^2\times\R^2$. \\

  \ni For the backward coefficients $g$ and $f$, we consider the following assumptions for the quadratic case.
  
\vspace{2mm}  
  
%
%
%
%

  

 \ni\textbf{(HBQ)} \begin{itemize}
 \item There exist three constants $M_g$, $K_g$ and $K_q$ such that 
  the functions $g$ and $f$ satisfy 
  \beqs
  |g(x)| & \leq & M_g\;,\\
  |g(x)-g(x')| & \leq & K_g |x-x'| \;,\\
  |f(t,x,y,z,u) - f(t,y',z,u)| & \leq & K_q |y-y'| \;,\\
  |f (t,x,y,z,u)|  & \leq & 
  K_q \big(1+ |y|+|z|^2+|u|\big) \;,
  \enqs
 for all $(t,x,x',y,y',z,u)\in [0,T] \times \R^2 \times \R^2 \times \R \times\R$. 
 \item For any $R >0$ there exists a function $mc^f_R$ such that $\lim_{\epsilon \rightarrow 0} mc^f_R(\epsilon) =0$ and
 \beqs
 |f(t,x,y,z, u-y) - f(t,x,y',z', u-y') | & \leq & mc^f_R(\epsilon)
 \enqs
 for all $(t,x,y,y',z,z',u) \in [0,T] \times\R \times  \R^2 \times \R^2 \times \R$ s.t. $|y|$, $|z|$, $|y'|$, $|z'| \leq R$ and $|y-y'| + |z-z'| \leq \epsilon$. 
 \item $f(t,.,u)=f(t,.,0)$ for all $u \in \R$ and all $t \in (\tau \wedge T, T]$.
 \item The function $f(t,x,y,.,u)$ is convexe (or concave) uniformly in  $(t,x,y,u)\in [0,T] \times \R
  \times\R\times\R$.  \end{itemize}
  
\vspace{2mm}  


\ni In the sequel $K$ denotes a generic constant appearing in \textbf{(HBL)}, \textbf{(HBQ)}  and \textbf{(HF)} and which may vary from line to line.\\

\ni In the purpose to prove the existence of a solution to the FBSDE (\ref{FSDE})-(\ref{BSDE}) we follow the decomposition approach initiated by \cite{kl11} and for that we introduce the recursive system of  FBSDEs associated with 
\reff{FSDE}-\reff{BSDE}.

\vspace{2mm}

\ni$\bullet$ Find $(X^1(\theta),Y^1(\theta),Z^1(\theta))\in\Sc^2_\F[0,T]\times\Sc_\F^\infty[\theta ,T]\times H_\F^2[\theta,T]$ such that 
\beq\label{FSDE1} 
  X^1_t(\theta)  & = & x + \int_0^t\hspace{-2mm} b\big(s, X^1_s(\theta)\big) ds + \int_0^t \hspace{-2mm}\sigma\big(s, X^1_s(\theta)\big) dW_s + \beta\big(\theta, X^1_{\theta^-} (\theta)\big) \1_{\theta \leq t} \;,~0\leq t\leq T\;,\\
\label{BSDE1}  Y^1_t\big(\theta\big)  & = & g\big(X^1_T(\theta)\big) + \int_t^T\hspace{-3mm} f\big(s, X^1_s(\theta), Y^1_s(\theta), Z^1_s(\theta), 0\big) ds - \int_t^T \hspace{-3mm}Z^1_s(\theta) dW_s \;,~\theta\leq t\leq T\;,~\quad\quad 
\enq
for all $\theta\in[0,T]$.

\vspace{3mm}

\ni$\bullet$ Find $(X^0,Y^0,Z^0)\in\Sc^2_\F[0,T]\times\Sc_\F^\infty[0 ,T]\times H_\F^2[0,T]$ such that 
\beq\label{FSDE0}
  X^0_t  & = & x + \int_0^t b(s, X^0_s) ds + \int_0^t \sigma(s, X^0_s) dW_s \;, \quad0\leq t\leq T\;,\\ \label{BSDE0}
  Y^0_t  & = & g(X^0_T) + \int_t^T f \big(s, X^0_s, Y^0_s, Z^0_s, Y^1_s(s) - Y^0_s \big) ds - \int_t^T Z^0_s dW_s \;, ~~0\leq t\leq T\;.\quad
\enq
Then, the link between the FBSDE \reff{FSDE}-\reff{BSDE} and the recursive system of FBSDEs \reff{FSDE1}-\reff{BSDE1} and \reff{FSDE0}-\reff{BSDE0} is given by the following result.
\begin{Theorem}  \label{theoreme existence unicite}
Assume that \textbf{(DH)}, \textbf{(HBI)}, \textbf{(H)}, \textbf{(HF)} and  \textbf{(HBL)} or  \textbf{(HBQ)} 
hold true. Then, the FBSDE \reff{FSDE}-\reff{BSDE} admits a unique solution $(X,Y,Z,U)\in\Sc^2_\G[0,T]\times\Sc^\infty_\G[0,T]\times H^2_\G[0,T]\times L^2(\lambda)$ given by
\begin{equation}\label{expression solution} \left\{
\begin{aligned}
X_t ~ & =  ~X^0_t \1_{t < \tau} + X^1_t(\tau) \1_{\tau \leq  t} \;,\\
Y_{t}  ~& =~ Y^0_t \mathds{1}_{t< \tau} + Y^1_t(\tau) \mathds{1}_{\tau \leq t} \;,\\
Z_{t} ~& =  ~Z^0_t \mathds{1}_{t \leq \tau} + Z^1_t(\tau) \mathds{1}_{\tau <  t} \;,\\
U_{t}~ & = ~ \big(Y^1_t(t) -Y^0_t\big) \mathds{1}_{t \leq \tau} \;,
\end{aligned}
\right.\end{equation}
 where  $(X^1(\theta),Y^1(\theta),Z^1(\theta))$ is the unique solution to the FBSDE \reff{FSDE1}-\reff{BSDE1} in $\Sc^2_\F[0,T]\times\Sc^\infty_\F[\theta,T]\times H^2_\F[\theta,T]$, for $\theta\in[0,T]$, and $(X^0,Y^0,Z^0)$ is the unique solution to the FBSDE \reff{FSDE0}-\reff{BSDE0} in $\Sc^2_\F[0,T]\times\Sc^\infty_\F[0,T]\times H^2_\F[0,T]$. 
\end{Theorem}
\ni\textbf{Proof.}

\ni\textbf{Step 1.} Solution to \reff{FSDE} under \textbf{(HF)}.

\ni Under  \textbf{(HF)} there exist unique processes $X^0\in\Sc^2_\F[0,T]$ satisfying \reff{FSDE0}, and  $X^1(\theta)\in\Sc^2_\F[0,T]$ satisfying \reff{FSDE1} for all $\theta\in [0,T]$ such that $X^1$ is $\Pc\Mc(\F)\otimes\Bc(\R_+)$-measurable. Then, from the definition of $H$,  we check that the process $X$ defined by 
\beq\label{decompX}
X_t & = & X^0_t\1_{t<\tau}+X^1_t(\tau)\1_{t\geq \tau}\;,\quad 
\enq
satisfies \reff{FSDE}. 
We now check that  $X\in\Sc^2_\G[0,T]$.
We first notice that from \textbf{(HF)}, there exists a constant $K$ such that 
 \beq\label{majS2X0}
 \E\Big[\sup_{t\in[0,T]}\big|X^0_t\big|^2 \Big]& \leq &  K\;. 
 \enq
 Then, from the definition of $X^0$ and $X^1$, we have for all $t \in [\theta, T]$
 \beqs
 \sup_{s\in[\theta,t]}\big|X^1_s(\theta)\big|^2 & \leq &  K\Big( \big| X^0_\theta \big|^2 +  \big|\beta(\theta,X^0_\theta)\big|^2+ \int_\theta^t \big| b(u,X^1_u(\theta))\big|^2du+  \sup_{s\in[\theta,t]} \Big|\int_\theta^s\sigma(u,X^1_u(\theta))dW_u\Big|^2\Big)\;. 
 \enqs
Using \textbf{(HF)} and BDG-inequality, we get 
 \beqs
\E\Big[ \sup_{s\in[\theta,t]}\big|X^1_s(\theta) \big|^2\Big] & \leq &  K\Big(1+  \int_\theta^t\E\Big[ \sup_{u\in[\theta,s]}\big|X^1_u(\theta) \big|^2\Big] du \Big)\;.
 \enqs
Applying Gronwall's lemma, we get 
 \beq\label{majS2X1}
 \sup_{\theta\in[0,T]}\big\|X^1(\theta)  \big\|_{\Sc_\F^2[\theta,T]} & \leq & K\;.
 \enq
 Combining \reff{decompX}, \reff{majS2X0} and \reff{majS2X1}, we get that $X\in \Sc_\G^2[0,T]$. 
Moreover still using \textbf{(HF)} we get the uniqueness of a solution to \reff{FSDE}   in $\Sc^2_\G[0,T]$. 

\vspace{2mm}

\ni\textbf{Step 2.} Solution to \reff{BSDE} under \textbf{(DH)}, \textbf{(HBI)}, \textbf{(H)} and \textbf{(HBL)}. 

\ni  To follow the decomposition approach initiated by the authors in \cite{kl11}, we need the generator to be predictable. To this end, we notice that in the BSDE \reff{BSDE}, we can replace the generator $(t,y,z,u)\mapsto f(t,X_t,y,z, (1-H_t)u)$ by the predictable map $(t,y,z,u)\mapsto f(t,X_{t^-},y,z,(1-H_{t^-})u)$. 

 Using the decomposition \reff{decompX}, we are able to write explicitly the  
decompositions of the $\Gc_T$-measurable random variable $g(X_T)$ and the $\Pc(\G)\otimes\Bc(\R)\otimes\Bc(\R)\otimes\Bc(\R)$-measurable map $(\omega,t,y,z,u)\mapsto f(t,X_{t^-}(\omega),y,z,u(1-H_{t^-}(\omega)))$ given by Lemma 2.1 in \cite{kl11}
\beqs
g(X_T) & = & g(X^0_T)   \1_{T<\tau}+ g(X^1_T(\tau))\1_{T\geq\tau}\;,\\
f(t,X_{t^-},y,z,(1-H_{t^-})u) & = & f^0(t,y,z,u) \1_{t\leq\tau}+ f^1(t,y,z,u,\tau)\1_{t>\tau} \;,
\enqs
with $f^0(t,y,z,u) = f(t,X^0_t,y,z,u)$ and  $f^1(t,y,z,u,\theta) = f(t,X^1_{t^-}(\theta),y,z,0)$, for all $(t,y,z,u,\theta)\in[0,T]\times\R\times\R\times\R\times\R_+$. 


\vspace{2mm}

Suppose now that \textbf{(DH)}, \textbf{(HBI)}, \textbf{(H)} and \textbf{(HBL)} hold true. Then, from Theorem C.1 in \cite{kl11}, the BSDE \reff{BSDE1} admits a $\Pc(\F)\otimes\Bc([0,T])$-measurable solution $(Y^1,Z^1)$ and the BSDE \reff{BSDE0}  admits a solution $(Y^0,Z^0)$. Using Proposition 2.1 in \cite{k00}, we obtain 
\beqs
{\|Y^1(\theta)\|}_{\Sc^\infty[\theta,T]} + {\|Z^1(\theta)\|}_{H^2[\theta,T]} & \leq  & K\;,
\enqs
for all $\theta\in[0,T]$, and 
\beqs
{\|Y^0\|}_{\Sc^\infty[0,T]} + {\|Z^0\|}_{H^2[0,T]} & \leq  & K\;.
\enqs
We can then apply Theorem 3.1 in \cite{kl11} and we get the existence of a solution to \reff{BSDE} in $\Sc^\infty_\G[0,T]\times H^2_\G[0,T]\times L^2(\lambda)$. 

Let $(Y,Z,U)$ and $(Y',Z',U')$ be two solutions to \reff{BSDE} in $\Sc^\infty_\G[0,T]\times H^2_\G[0,T]\times L^2(\lambda)$. Since 
$f(t,x,y,z,(1-H_t)u)=f(t,x,y,z,0)$ for all $t\in(\tau\wedge T,T]$ and $\lambda$ vanishes on $(\tau\wedge T,T]$, we can assume w.l.o.g. that $U_t=U_t'=0$ for $t\in(\tau\wedge T,T]$. Then, from \textbf{(DH)}, \textbf{(HBI)}, \textbf{(H)} and \textbf{(HBL)}, we can apply Theorem 4.1 in \cite{kl11} and we get that $Y\leq Y'$. Since $Y$ and $Y'$ play the same role, we obtain $Y=Y'$. Identifying the pure jump parts of $Y$ and $Y'$  gives $U=U'$. Finally, identifying the unbounded variation gives $Z=Z'$.

\vspace{2mm}

\vspace{2mm}

\ni\textbf{Case 3.}  Solution to \reff{BSDE} under \textbf{(DH)}, \textbf{(HBI)}, \textbf{(H)} and  \textbf{(HBQ)}.

\ni The existence of a solution $(Y,Z,U)\in\Sc^2_\G[0,T]\times L^2_\G[0,T]\times L^2(\lambda)$ is a direct consequence of Proposition 3.1 in \cite{kl11}.
We then notice that from the definition of $H$ we have $f(t,x,y,z,u(1-H_t))=f(t,x,y,z,0)$ for all $t\in(\tau\wedge T,T]$.  This property and \textbf{(DH)}, \textbf{(HBI)}, \textbf{(H)} and \textbf{(HBQ)}  allow to apply  Theorem 4.2 in \cite{kl11}, which gives the uniqueness of a solution of \reff{BSDE}.   
\ep

\vspace{2mm}

Throughout the sequel, we give an approximation of the solution to the FBSDE \reff{FSDE}-\reff{BSDE} by studying the approximation of the solutions to the recursive system of FBSDEs \reff{FSDE1}-\reff{BSDE1} and \reff{FSDE0}-\reff{BSDE0}. 

\section{Discrete-time scheme for the FBSDE}
\setcounter{equation}{0} \setcounter{Assumption}{0}
\setcounter{Theorem}{0} \setcounter{Proposition}{0}
\setcounter{Corollary}{0} \setcounter{Lemma}{0}
\setcounter{Definition}{0} \setcounter{Remark}{0}

In this section, we introduce a discrete-time approximation of the solution $(X,Y,Z,U)$ to the FBSDE \reff{FSDE}-\reff{BSDE} based on its decomposition given by Theorem \ref{theoreme existence unicite}. 

Throughout the sequel, we consider a discretization grid $\pi:=\{t_0,\ldots,t_n\}$ of $[0,T]$ with $0=t_0<t_1<\ldots<t_n=T$.  For $t\in[0,T]$, we denote by $\pi(t)$ the largest element of $\pi$ smaller than $t$
\beqs
\pi(t) & := & \max\big\{\;t_i\;,~i=0,\ldots,n~|~t_i\leq t\;\big\}\;.
\enqs
We also denote by $|\pi|$ the mesh of $\pi$
\beqs
|\pi| & := & \max\big\{\;t_{i+1}-t_i\;,~i=0,\ldots,n-1\;\big\}\;,
\enqs
that we suppose satisfying $|\pi|\leq 1$, 
and by $\Delta W^\pi_i$ (resp. $\Delta t^\pi_i $) the increment of $W$ (resp. the difference) between $t_i$ and $t_{i-1}$: $\Delta W^\pi_i := W_{t_{i}} - W_{t_{i-1}}$ (resp. $\Delta t^\pi_i :=t_i-t_{i-1}$), 
for $1\leq i\leq n$.

\subsection{Discrete-time scheme for $X$ }

\ni We introduce an approximation of the process $X$ based on the discretization of the processes $X^0$ and $X^1$.

\vspace{2mm}

\ni$\bullet$ \emph{Euler scheme for $X^0$}. We consider the scheme $X^{0,\pi}$ defined by
\begin{equation*} \left\{
\begin{aligned}
  \XzePi_{t_0} ~ & = ~ x \;,\\
  \XzePi_{t_i} ~ & = ~   \XzePi_{t_{i-1}} + b(t_{i-1}, \XzePi_{t_{i-1}})\Delta t_i^\pi + \sigma(t_{i-1}, \XzePi_{t_{i-1}}) \Delta W^\pi_i  \;,\quad 1\leq i\leq n\;.\\
   \end{aligned}
\right.\end{equation*}
\ni$\bullet$ \emph{Euler scheme for $X^1$}. Since the process $X^1$ depends on two parameters $t$ and $\theta$, we introduce a discretization of $X^1$ in these two variables. We then consider the following scheme
\begin{equation}\label{Euler1} \left\{
\begin{aligned}
  \XunPi_{t_0}(\pi(\theta)) ~ & = ~ x+\beta(t_0,x)\1_{\pi(\theta) = 0} \;,\\
  \XunPi_{t_i}(\pi(\theta)) ~ & = ~   \XunPi_{t_{i-1}}(\pi(\theta)) + b(t_{i-1}, \XunPi_{t_{i-1}}(\pi(\theta)))\Delta t^\pi_i + \sigma( t_{i-1},\XunPi_{t_{i-1}}(\pi(\theta))) \Delta W^\pi_i\\
   & \quad  + \beta(t_{i-1}, \XunPi_{t_{i-1}}(\pi(\theta))) \1_{t_i = \pi(\theta)} \;,\quad\ 1\leq i\leq n \;, \quad 0 \leq \theta \leq T\;.\\
   \end{aligned}
\right.\end{equation}

\ni We are now able to provide an approximation of the process $X$ solution to the FSDE \reff{FSDE}.
We consider the scheme $X^\pi$ defined by 
\beq\label{EulerX}
X^\pi_{t} & = & \XzePi_{\pi(t)}\1_{t < \tau} +\XunPi_{\pi(t)}(\pi(\tau))\1_{t \geq \tau}\;,\quad 0 \leq t \leq T\;.
\enq
We shall denote by $\{\Fc^{0,\pi}_i\}_{0 \leq i \leq n}$ (resp. $\{\Fc^{1,\pi}_i(\theta)\}_{0 \leq i \leq n}$) the discrete-time filtration  associated with $X^{0,\pi}$ (resp. $X^{1,\pi}$)  
\beqs
\Fc^{0,\pi}_i& := & \sigma( \XzePi_{t_j},\; j \leq i) \\
 \mbox{ (resp. }\Fc^{1,\pi}_i(\theta)  & := &  \sigma( \XunPi_{t_j}(\theta),\; j \leq i)\mbox{)}\;.%
\enqs

\subsection{Discrete-time scheme for $(Y, Z, U)$ }

We introduce an approximation of $(Y,Z)$ based on the discretization of $(Y^0,Z^0)$ and $(Y^1,Z^1)$.   
To this end we introduce the backward implicit schemes on $\pi$ associated with the BSDEs \reff{BSDE1} and \reff{BSDE0}. Since the system is recursively coupled, we first introduce the scheme associated with  \reff{BSDE1}. We then use it to define the scheme associated with \reff{BSDE0}. 

\vspace{2mm}

\ni$\bullet$   \textit{Backward Euler scheme for $(Y^1,Z^1)$}. We consider the implicit scheme $(\YunPi,\ZunPi)$ defined by 
\begin{equation}\label{YZ1} \left\{
\begin{aligned}
  \YunPi_{T}(\pi ( \theta )) ~ & = ~ g ( \XunPi_{T}(\pi ( \theta )) ) \;,\\
  \YunPi_{t_{i-1}}(\pi ( \theta )) ~ & = ~ \E^{1,\pi(\theta)}_{i-1} \big[ \YunPi_{t_{i}}(\pi ( \theta )) \big] + f \big( t_{i-1},\XunPi_{t_{i-1}}(\pi ( \theta )), \YunPi_{t_{i-1}}(\pi ( \theta )), \ZunPi_{t_{i-1}}(\pi ( \theta )), 0\big) \Delta t_i^\pi \;,\\
    \ZunPi_{t_{i-1}}(\pi ( \theta )) ~ & = ~ \frac{1}{\Delta t_i^\pi } \E^{1,\pi(\theta)}_{i-1} \big[ \YunPi_{t_i}(\pi ( \theta )) \Delta W^\pi_{i} \big] \;,\quad  \pi(\theta) \leq t_{i-1}\;,~1\leq i\leq n\;,
    \end{aligned}
\right.\end{equation}
  where $\E^{1,s}_{i} = \E [ \; . \; | \Fc^{1,\pi}_{i}(s)]$ for $0\leq i\leq n$ and $s\in[0,T]$.

\vspace{2mm}

\ni$\bullet$  \textit{Backward Euler scheme for $(Y^0,Z^0)$}. Since the generator of \reff{BSDE0} involves the process ${(Y^1_t(t))}_{t\in[0,T]}$, we consider a discretization based on $\YunPi$. We therefore consider the scheme $(\YzePi,\ZzePi)$ defined by 
\begin{equation}\label{YZ0} \left\{
\begin{aligned}
  \YzePi_{T} ~ & = ~ g ( \XzePi_{T} ) \;,\\
  \YzePi_{t_{i-1}} ~ & = ~ \E^0_{i-1} \big[ \YzePi_{t_{i}} \big] + \bar f^\pi \big( t_{i-1},\XzePi_{t_{i-1}}, \YzePi_{t_{i-1}}, \ZzePi_{t_{i-1}}\big) \Delta t^\pi_{i} \;,\\
    \ZzePi_{t_{i-1}} ~ & = ~ \frac{1}{\Delta t^\pi_{i} } \E^0_{i-1} \big[ \YzePi_{t_i} \Delta W^\pi_{i} \big] \;,\quad1\leq i\leq n\;,
    \end{aligned}
\right.\end{equation}
  where $\E^0_{i} = \E [ \; . \; | \Fc^{0,\pi}_{i}]$ for $0\leq i\leq n$, and $\bar f^\pi$ is defined by
  \beqs
  \bar f^\pi(t,x,y,z) & = & f\big(t,x,y,z,\YunPi_{\pi(t)}(\pi(t))-y\big)\;,
  \enqs
  for all $(t,x,y,z)\in[0,T]\times\R \times\R\times\R$.
  
\ni   We then consider the following scheme for the solution $(Y,Z,U)$ of the BSDE \reff{BSDE}
\begin{equation}\label{SchemeYZ} \left\{
\begin{aligned}
Y^\pi_{t} ~& =~  Y^{0, \pi}_{\pi(t)} \mathds{1}_{t < \tau} + \YunPi_{\pi(t)} ( \pi ( \tau ) ) \mathds{1}_{t \geq  \tau} \;,\\
Z^\pi_{t}~ & =~  Z^{0, \pi}_{\pi(t)} \mathds{1}_{t \leq  \tau} + \ZunPi_{\pi(t)} ( \pi (  \tau ) ) \mathds{1}_{t >  \tau} \;,\\
U^\pi_{t} ~& =~  \big( \YunPi_{\pi(t)} ( \pi(t)) - Y^{0, \pi}_{\pi(t)} \big) \mathds{1}_{t \leq  \tau} \;,
    \end{aligned}
\right.\end{equation}
for $t\in[0,T]$.

\section{Convergence of the forward scheme}
\setcounter{equation}{0} \setcounter{Assumption}{0}
\setcounter{Theorem}{0} \setcounter{Proposition}{0}
\setcounter{Corollary}{0} \setcounter{Lemma}{0}
\setcounter{Definition}{0} \setcounter{Remark}{0}

We introduce the following assumption, which will be used to control the error between $X$ and $X^\pi$.

\vspace{2mm}

\ni\textbf{(HFD)} There exists a constant $K$ such that the functions $b$, $\sigma$ and $\beta$ satisfy
\beqs
\big|b(t,x)-b(t',x)\big| +\big|\sigma(t,x)-\sigma(t',x)\big| & \leq & K|t-t'|^{1\over 2}\;,\\
\big|\beta(t,x)-\beta(t',x)\big| +\big|\sigma(t,x)-\sigma(t',x)\big| & \leq & K|t-t'|\;,
\enqs
for all $(t,t',x)\in[0,T]\times[0,T]\times\R$.

\vspace{2mm}

In the following we provide an error estimate of the approximation schemes for $X^0$ and $X^1$ which are used to control the error between $X$ and $X^\pi$.




\subsection{Error estimates for $X^0$ and $X^1$}

\ni Under \textbf{(HF)} and \textbf{(HFD)}, the upper bound of the error between $X^0$ and its Euler scheme $X^{0, \pi}$ is well understood, see e.g. \cite{klopla92}, and we have
\beq\label{majerrX0}
\E\Big[ \sup_{t\in[0,T]} \big|X^0_t - X^{0,\pi}_{\pi(t)} \big|^2 \Big] & \leq & K | \pi | \;,
\enq
for some constant $K$ which does not depend on $\pi$.

\vspace{2mm}

\ni The next result provides an upper bound for the error between $X^1$ and its Euler scheme $\XunPi$ defined by \reff{Euler1}.  

\begin{Theorem}\label{erreur x1}
Under \textbf{(HF)} and \textbf{(HFD)}, we have the following  estimate 
\beqs
\sup_{\theta\in [0, T]} \E \Big[\sup_{t\in [\theta, T]} \big| X^1_t(\theta) - \XunPi_{\pi(t)}(\pi(\theta)) \big|^2 \Big]  & \leq &  K |\pi|  \;,
\enqs
for a constant $K$ which does not depend on $\pi$.
\end{Theorem}

\ni\textbf{Proof.}
Fix $\theta\in[0,T]$, we then have
\beq\nonumber
\E \Big[ \sup_{t\in [\theta, T]}\big| X^1_t(\theta) - \XunPi_{\pi(t)}(\pi(\theta)) \big|^2 \Big] &  \leq &  2\; \E \Big[ \sup_{t\in [\theta, T]}\big| X^1_t(\theta) - X^1_t(\pi(\theta)) \big|^2 \Big] \\\label{decompX-Xpi}
 & & +\; 2 \; \E \Big[ \sup_{t\in [\theta, T]}\big| X^1_t(\pi(\theta)) - \XunPi_{\pi(t)}(\pi(\theta)) \big|^2 \Big] \;.
\enq
We study separately the two terms of the right hand side.

\ni Since $\pi(\theta) \leq \theta \leq t$, we have by definition $X^1_{s} ( \pi ( \theta ) )  =  X^0_s$ for all $s\in [ 0 ,  \pi ( \theta))$,  and $ X^1_{s} ( \theta )=X^0_s$ for all $s\in [ 0 ,  \theta)$, which implies
\beqs
X^1_t ( \theta ) - X^1_t ( \pi ( \theta ) ) & = & \int_{\pi(\theta)}^{\theta} b \big(s, X^0_s\big) ds +  \int_{\pi(\theta)}^{\theta} \sigma\big(s, X^0_s\big) dW_s + \beta \big( \theta, X^0_\theta\big) +  \int_\theta ^t b\big(s,X^1_s ( \theta ) \big)ds\\
 & & + \int_\theta^t \sigma \big(s, X^1_s ( \theta ) \big) dW_s - \beta \big( \pi(\theta), X^0_{\pi(\theta)}\big) - \int_{\pi(\theta)}^{t} b\big(s, X^1_s ( \pi(\theta))\big) ds \\
 & &-  \int_{\pi(\theta)}^{t} \sigma\big(s, X^1_s ( \pi(\theta))\big) dW_s   \;,
\enqs
for all $t\in[\theta,T]$. 

\ni Hence, there exists a constant $K$ such that 
\beq
\nonumber \big| X^1_t ( \theta ) - X^1_t ( \pi ( \theta ) ) \big|^2 & \leq & K \Big \{ \Big|  \int_{\pi(\theta)}^{\theta} b\big(s, X^0_s\big) ds \Big|^2 +  \Big|  \int_{\pi(\theta)}^{\theta} b\big(s, X^1_s ( \pi(\theta))\big) ds \Big|^2  \\
&& \nonumber+  \int_{\theta}^t \Big| b\big(s,X^1_s ( \theta ) \big) - b \big(s, X^1_s ( \pi ( \theta ) ) \big) \Big|^2 ds  \\
& & \nonumber + \; \Big| \int_{\pi(\theta)}^{\theta} \sigma\big(s, X^0_s\big) dW_s \Big|^2 +  \Big| \int_{\pi(\theta)}^{\theta} \sigma\big(s, X^1_s ( \pi(\theta))\big) dW_s \Big|^2 \\
&& \nonumber + \; \Big| \int_\theta^t \Big( \sigma \big(s, X^1_s ( \theta ) \big) - \sigma \big(  s, X^1_s ( \pi ( \theta )  ) \big) \Big) dW_s \Big|^2 \\
& & + \; \big| \beta \big(\theta, X^0_{\theta}  \big) - \beta \big(\pi ( \theta ), X^0_{\pi ( \theta )}  \big) \big|^2 \Big\}  \label{deltaX1}   \;.
\enq
From \textbf{(HF)} and \textbf{(HFD)}, we have
\beqs
\E\big| \beta \big(\theta, X^0_{\theta}  \big) - \beta \big(\pi ( \theta ), X^0_{\pi ( \theta )}  \big) \big|^2 & \leq & K\big( |\pi|^2+\E \big|X^0_\theta-X^0_{\pi(\theta)}\big|^2\big)\;.
\enqs
We have from \textbf{(HF)} and \reff{majS2X0}
\beqs
\E \Big [\Big|  \int_{\pi(\theta)}^{\theta} b(s, X^0_s) ds \Big|^2 +  \Big| \int_{\pi(\theta)}^{\theta} \sigma(s, X^0_s) dW_s \Big|^2 \Big ] & \leq & K|\pi| \;,
\enqs
which implies in particular $\E |X^0_\theta-X^0_{\pi(\theta)}|^2 \leq K | \pi |$ and hence 
\beqs
\E |\beta(\theta,X^0_\theta)-\beta(\pi(\theta),X^0_{\pi(\theta)})|^2 \leq K | \pi |\;.
\enqs
 We have also from \textbf{(HF)} and \reff{majS2X1}
\beqs
\E \Big [ \Big|  \int_{\pi(\theta)}^{\theta} b\big(s, X^1_s ( \pi(\theta))\big) ds \Big|^2  +  \Big| \int_{\pi(\theta)}^{\theta} \sigma\big(s, X^1_s ( \pi(\theta))\big) dW_s \Big|^2 \Big] & \leq & K|\pi| \;.
\enqs
Combining these inequalities with \reff{deltaX1}, \textbf{(HF)} and BDG-inequality, we get
\beqs
\E \Big[ \sup_{u\in[\theta,t]}\big| X^1_u ( \theta ) - X^1_{ u }( \pi(  \theta ) ) \big|^2 \Big] & \leq & K \Big(\int_\theta ^t \E \Big[ \sup_{u\in[\theta,s]}\big| X^1_u ( \theta ) - X^1_{ u }( \pi(  \theta ) ) \big|^2 \Big] ds + |\pi|\Big) \;.
\enqs
Applying Gronwall's lemma, we get
\beq\label{maj1DX1}
\E \Big[ \sup_{t\in[\theta,T]}\big| X^1_t ( \theta ) - X^1_{ t }( \pi(  \theta ) ) \big|^2 \Big] & \leq & K |\pi|  \;.
\enq

\vspace{2mm}

\ni To find an upper bound for the term $\E [ \sup_{t \in [ \theta, T ]} | X^1_t(\pi(\theta)) - \XunPi_{\pi(t)}(\pi(\theta)) |^2 ]$ we introduce the scheme $\XtilPi_.( \pi ( \theta ) )$ defined by
\begin{equation*} \left\{
\begin{aligned}
  \XtilPi_{\pi ( \theta ) }(\pi ( \theta )) ~ & = ~ X^1_{\pi ( \theta ) }(\pi ( \theta )) \;,\\
   \XtilPi_{t_i}(\pi ( \theta )) ~ & = ~ \XtilPi_{t_{i-1}}(\pi ( \theta )) + b\big( t_{i-1}, \XtilPi_{t_{i-1}}(\pi ( \theta ))\big)\Delta t^\pi_i + \sigma\big( t_{i-1},\XtilPi_{t_{i-1}} (\pi ( \theta ))\big) \Delta W^\pi_i \;, ~ t_i > \pi ( \theta ) \;.\\
   \end{aligned}
\right.\end{equation*}
We have the inequality
\beq\label{decomp2DX2}
\hspace{-7mm} \E \Big[ \sup_{t\in[\theta,T]}\big| X^1_t(\pi ( \theta )) - \XunPi_{\pi(t)}( \pi ( \theta )) \big|^2 \Big] & \leq  & 2 \; \E \Big[ \sup_{t\in[\theta,T]}\big| X^1_t( \pi (\theta)) - \XtilPi_{\pi(t)} (\pi(\theta)) \big|^2 \Big] \nonumber\\
& & + \; 2 \; \E \Big[ \sup_{t\in[\theta,T]}\big| \XtilPi_{\pi(t)} (\pi(\theta)) - \XunPi_{\pi(t)}(\pi(\theta)) \big|^2 \Big] \;.
\enq
Since $\XtilPi( \pi ( \theta ) )$ is the Euler scheme of $X^1( \pi ( \theta ) )$ on $[\pi(\theta),T]$,  we have under \textbf{(HF)} and \textbf{(HFD)} (see e.g. \cite{klopla92})
\beqs
\E \Big[ \sup_{t\in[\theta,T]}\big| X^1_t ( \pi (\theta)) - \XtilPi_{\pi(t)} (\pi(\theta)) \big|^2 \Big] & \leq & K\Big(1+\E\Big[\big|X^1_{\pi(\theta)}(\pi(\theta))\big|^2\Big]\Big) |\pi| \;,
\enqs
for some constant $K$ which neither depends on $\pi$ nor on $\theta$. From \reff{majS2X1}, we get 
\beq\label{majDX21}
\E \Big[ \sup_{t\in[\theta,T]}\big| X^1_t( \pi (\theta)) - \XtilPi_{\pi(t)} (\pi(\theta)) \big|^2 \Big] & \leq & K
|\pi|\;,
\enq
for all $\theta\in[0,T]$.

We now study the term $ \E \big[ \sup_{t\in[\theta,T]}\big| \XtilPi_{\pi(t)} (\pi(\theta)) - \XunPi_{\pi(t)}(\pi(\theta)) \big|^2 \big]$. We first notice that we have the following identity
\beqs
 \E \Big[ \sup_{t\in[\theta,T]}\big| \XtilPi_{\pi(t)} (\pi(\theta)) - \XunPi_{\pi(t)}(\pi(\theta)) \big|^2 \Big] & = & \E \Big[ \sup_{t\in[\pi(\theta),T]}\big| \XtilPi_{\pi(t)} (\pi(\theta)) - \XunPi_{\pi(t)}(\pi(\theta)) \big|^2 \Big]\;.
 \enqs
 Hence we can work with the second term. 
  From the definition of $\XtilPi$ and $\XunPi$, we get
\begin{multline*}
\sup_{u\in[\pi(\theta),t]}\big| \XtilPi_{\pi(u)} (\pi(\theta)) - \XunPi_{\pi(u)}(\pi(\theta)) \big|^2  \leq \\
K\Big(\big| X^1_{\pi(\theta)} (\pi(\theta)) - \XunPi_{\pi(\theta)}(\pi(\theta))\big|^2 
+\int_{\pi(\theta)}^{\pi(t)}\Big|b\big(\pi(s),\XtilPi_{\pi(s)}(\pi(\theta))\big)-b\big(\pi(s),\XunPi_{\pi(s)}(\pi(\theta))\big)\Big|^2ds   \\
+\sup_{u\in[\pi(\theta),t]}\Big|\int_{\pi(\theta)}^{\pi(u)}\hspace{-2mm}\Big(\sigma\big(\pi(s),\XtilPi_{\pi(s)}(\pi(\theta))\big)-\sigma\big(\pi(s),\XunPi_{\pi(s)}(\pi(\theta))\big)\Big)dW_s\Big|^2\Big)\;. 
\end{multline*}
Then, using \textbf{(HF)} and BDG-inequality, we get 
\beqs
\E\Big[\sup_{u\in[\pi(\theta),t]}\big| \XtilPi_{\pi(u)} (\pi(\theta)) - \XunPi_{\pi(u)}(\pi(\theta)) \big|^2\Big]\hspace{-3mm} & \leq & \hspace{-3mm} K\Big(\E\big| X^1_{\pi(\theta)} (\pi(\theta)) - \XunPi_{\pi(\theta)}(\pi(\theta))\big|^2 
 \\
 & & \hspace{-3mm}  +  \hspace{-1mm}  \int_{\pi(\theta)}^{t}  \hspace{-2mm} \E\Big[ \sup_{u\in[\pi(\theta),s]}\big| \XtilPi_{\pi(u)} (\pi(\theta)) - \XunPi_{\pi(u)}(\pi(\theta)) \big|^2\Big]ds\Big) \,.
\enqs
From Lipschitz property of $\beta$, we have 
\beqs
\E\big| X^1_{\pi(\theta)} (\pi(\theta)) - \XunPi_{\pi(\theta)}(\pi(\theta))\big|^2 & = & \E\big| X^0_{\pi(\theta)} +\beta\big(\pi(\theta),X^0_{\pi(\theta)} \big) - \XzePi_{\pi(\theta)}-\beta\big(\pi(\theta),\XzePi_{\pi(\theta)}\big)\big|^2\\
 & \leq & K\,\E\big| X^0_{\pi(\theta)} - \XzePi_{\pi(\theta)}\big|^2\;.
\enqs
This last inequality with \reff{majerrX0} gives
\beqs
\E\big| X^1_{\pi(\theta)} (\pi(\theta)) - \XunPi_{\pi(\theta)}(\pi(\theta))\big|^2 & \leq & K|\pi|\;.
\enqs
Applying Gronwall's lemma, we get 
\beq\label{majDX22}
\E\Big[\sup_{t\in[\pi(\theta),T]}\big| \XtilPi_{\pi(t)} (\pi(\theta)) - \XunPi_{\pi(t)}(\pi(\theta)) \big|^2\Big] & \leq & K|\pi|\;.
\enq
Combining \reff{decompX-Xpi}, \reff{maj1DX1}, \reff{decomp2DX2}, \reff{majDX21} and \reff{majDX22}, we get the result.\ep


\subsection{Error estimate for the FSDE with a jump}

We are now able to provide an estimate of the error approximation of the process $X$ by its scheme $X^\pi$ defined by \reff{EulerX}. 
\begin{Theorem}
Under \textbf{(HF)} and \textbf{(HFD)}, we have the following  estimate 
\beqs
\E \Big[ \sup_{t\in[0,T]} \big| X_{t} - X^\pi_{t} \big|^2 \Big] & \leq &  K |\pi|  \;,
\enqs
for a constant $K$ which does not depend on $\pi$.
\end{Theorem}
\ni\textbf{Proof.} From the definition of $X^\pi$, \textbf{(DH)} and \reff{majerrX0} we have 
\beqs
\E\Big[  \sup_{t\in[0,T]}\big| X_{t} - X^\pi_{t} \big|^2 \Big] & \leq & \hspace{-2mm} \E\Big[  \sup_{t\in[0,\tau)} \big| X^0_{t} - X^{0,\pi}_{\pi(t)} \big|^2 \Big] +  \E\Big[  \sup_{t\in[\tau,T]}\big| X^1_{t}(\tau) - X^{1,\pi}_{\pi(t)}(\pi(\tau)) \big|^2 \Big]\\
 & \leq  &  \hspace{-2mm} \E\Big[  \sup_{t\in[0,T]}  \hspace{-1mm} \big| X^0_{t} - X^{0,\pi}_{\pi(t)} \big|^2\Big] + \hspace{-1mm}  \int_0^T \hspace{-2mm}\E\Big[  \sup_{t\in[\theta,T]} \hspace{-1mm}\big| X^1_{t}(\theta) - X^{1,\pi}_{\pi(t)}(\pi(\theta)) \big|^2\gamma_{T}(\theta) \Big]d\theta \\
 & \leq & K\Big(|\pi|+ \sup_{\theta\in[0,T]}\E\Big[ \sup_{s\in[\theta,T]}\big| X^1_{s}(\theta) - X^{1,\pi}_{\pi(s)}(\pi(\theta)) \big|^2\Big]\Big) \;.
\enqs 
From Theorem \ref{erreur x1}, we get
\beqs
\E\Big[ \sup_{t\in[0,T]}\big| X_{t} - X^\pi_{t} \big|^2 \Big] & \leq & K|\pi|\;.
\enqs
\ep
\section{Convergence of the backward scheme in the Lipschitz case}\label{lipschitz}

\setcounter{equation}{0} \setcounter{Assumption}{0}
\setcounter{Theorem}{0} \setcounter{Proposition}{0}
\setcounter{Corollary}{0} \setcounter{Lemma}{0}
\setcounter{Definition}{0} \setcounter{Remark}{0}
To provide  error estimates for the Euler scheme of the BSDE, we need an additional regularity property for the coefficients $g$ and $f$. We then introduce the following assumption. 

\vspace{2mm}

\ni\textbf{(HBLD)} There exists a constant $K$ such that the functions $g$ and $f$ satisfy
\beqs
\big|g(x)-g(x')\big|+ \big|f(t,x,y,z,u)-f(t',x',y,z,u)\big| & \leq & K\big(|x-x'|+|t-t'|^{1\over 2}\big)\;,
\enqs
for all $(t,t',x,x',y,z,u)\in[0,T]^2\times\R^2\times\R\times\R\times\R$.

\vspace{2mm}

\ni We are now ready to provide error estimates of the approximation schemes for $(Y^0,Z^0)$ and $(Y^1,Z^1)$, and then for $(Y,Z)$.

\subsection{Regularity results}

In this part, we give some results on the regularity of the processes $Z^1$ and $Z^0$. We denote $\Fc^0_t := \sigma\{ X^0_s\;, ~0 \leq s \leq t \}$ and $\Fc^1_t(\theta) := \sigma\{ X^1_s(\theta)\;, ~\theta \leq s \leq t \}$. 

\begin{Proposition}\label{regularite z1}
Under \textbf{(HF)}, \textbf{(HFD)}, \textbf{(HBL)} and \textbf{(HBLD)}, there exists a constant $K$ such that 
\beq\label{borne-reg-Z1}
\E \Big[ \int_\theta^T \big| Z^1_t(\theta) - Z^1_{\pi(t)}(\theta) \big|^2 dt \Big] & \leq & K \Big( 1 + \E \Big[ \big| X^1_\theta ( \theta) \big|^4\Big]^{1\over2}\Big) |\pi| \;,
\enq
for all $\theta \in \pi$.
\end{Proposition}

\ni\textbf{Proof.} 
We first suppose that $b$, $\sigma$, $f$ and $g$ are in $C^1_b$. Let us define the processes $\Lambda$ and $M$ by
\beqs
\Lambda_t &: = & \exp \Big( \int_\theta^t \partial_y f \big( \Theta^1_r ( \theta) \big) dr \Big)\;,
\enqs
and
\beqs
M_t & := & 1 + \int_\theta^t M_r \partial_z  f \big( \Theta^1_r ( \theta) \big) dW_r \;,
\enqs
where $\Theta^1_r ( \theta) := (r, X^1_r ( \theta),Y^1_r ( \theta),Z^1_r ( \theta) , 0)$. We give classically the link between $\nabla^\theta X^1_t(\theta) (:= \partial X^1_t(\theta) / \partial X^1_\theta ( \theta))$ and $(D_s X^1_t(\theta))_{\theta \leq s \leq t}$ the Malliavin derivative of $X^1_t(\theta)$. Recall that $X^1(\theta)$ satisfies
\beqs
X^1_t(\theta)  & = & X^1_\theta(\theta)  + \int_\theta^t\hspace{-2mm} b\big(r, X^1_r(\theta)\big) dr + \int_\theta^t \hspace{-2mm}\sigma\big(r, X^1_r(\theta)\big) dW_r \;,\quad \theta\leq t \leq T\;.
\enqs
Therefore, we get
\beqs
\nabla^\theta X^1_t(\theta) & = & 1 + \int_\theta^t \partial_x b \big(r, X^1_r(\theta)\big) \nabla^\theta X^1_r(\theta) dr + \int_\theta^t \partial_x \sigma \big(r, X^1_r(\theta)\big)  \nabla^\theta X^1_r(\theta) dW_r \;,\quad \theta\leq t \leq T\;,
\enqs
 and for $\theta \leq s \leq t$
\beqs
D_s  X^1_t(\theta) & = & \sigma\big( s, X^1_s(\theta)\big) + \int_s^t \partial_x b \big( r, X^1_r(\theta)\big) D_s  X^1_r(\theta) dr + \int_s^t \partial_x \sigma\big( r, X^1_r(\theta)\big) D_s  X^1_r(\theta) dW_r \;.
\enqs
Thus, we have
\beq \label{lien D et nabla}
D_s  X^1_t(\theta) & = & \nabla^\theta X^1_t(\theta)  \big[ \nabla^\theta X^1_s(\theta) \big]^{-1}  \sigma\big( s, X^1_s(\theta)\big) \;.
\enq
Using Malliavin calculus we obtain that a version of $Z^1(\theta)$ is given by ${(D_{t}Y_t^1(\theta))}_{t\in[\theta,T]}$. 
By It\^{o}'s formula, we get 
\beqs
\Lambda_t M_t Z^1_t( \theta ) & = & \E \Big[ M_T \Big( \Lambda_T \nabla g \big(X^1_T ( \theta)\big) D_t X^1_T ( \theta) + \int_t^T \partial_x f \big( \Theta^1_r ( \theta) \big) D_t X^1_r ( \theta) \Lambda_r dr \Big)\Big| \Fc^1_t (\theta) \Big] \;,
\enqs
for $t \in [\theta, T]$. Using \reff{lien D et nabla}, we get
\beqs
\Lambda_t M_t Z^1_t( \theta ) & = & \E \Big[ M_T \Big( \Lambda_T \nabla g \big(X^1_T ( \theta)\big) \nabla^\theta X^1_T ( \theta) + \int_t^T F_r \Lambda_r dr\Big) \Big| \Fc^1_t(\theta) \Big]  \big[ \nabla^\theta X^1_t(\theta) \big]^{-1}  \sigma\big( t, X^1_t(\theta)\big) \;,
\enqs
with $F_r := \partial_x f \big( \Theta^1_r ( \theta) \big) \nabla^\theta X^1_r ( \theta)$. This implies that
\beqs
\Lambda_t M_t Z^1_t( \theta ) & = & \Big( \E [G | \Fc^1_t(\theta)]  - M_t\int_\theta^t F_r \Lambda_r dr \Big)  \big[ \nabla^\theta X^1_t(\theta) \big]^{-1}  \sigma\big( t, X^1_t(\theta)\big) \;,
\enqs
with $G := M_T \Big(\Lambda_T \nabla g (X^1_T ( \theta)) \nabla^\theta X^1_T ( \theta) + \int_\theta^T F_r \Lambda_r dr\Big)$.
Since $b$, $\sigma$, $f$ and $g$ have bounded derivatives, we have 
\beq\label{momG}
\E\big[|G|^p\big] & < & \infty\;, \quad p\geq 2\;. 
\enq
Define $m_r := \E [ G | \Fc^1_r(\theta) ]$ for $r\in[\theta,T]$. From \reff{momG}  and  Doob's inequality, we have
\beq\label{momm}
\| m \|_{S^p[\theta,T]} & < & \infty\;, \quad p\geq 2\;. 
\enq
Hence, there exists a process $\phi$ such that 
\beqs
m_r &= & \E[G | \Fc^1_\theta(\theta) ] + \int_\theta^r \phi_u dW_u\;,\quad r\in[\theta,T] \;,
\enqs
and 
\beqs
\| \phi \|_{H^p[\theta,T]} & < & \infty\;, \quad p\geq 2\;. 
\enqs
We define $\tilde Z$ by
\beqs
\tilde Z_t(\theta) & := & (\Lambda_t M_t )^{-1} \Big( m_t - M_t \int_\theta^t F_r \Lambda_r dr \Big) \big[ \nabla^\theta X^1_t ( \theta ) \big]^{-1} \;.
\enqs
By It\^{o}'s formula, we can write
\beqs
\tilde Z _t (\theta)& = & \tilde Z_\theta(\theta) + \int_\theta^t \alpha^1_r(\theta) dr + \int_\theta^t \alpha^2_r(\theta) dW_r\;,\quad \theta\leq r\leq T\;.
\enqs
Since $b$, $\sigma$, $f$ and $g$ have bounded derivatives, we get from \reff{momm} 
\beq\label{avantthomasse}
\sup_{\theta\in[0,T]}\| \tilde Z(\theta) \|^p_{\Sc^p [\theta, T ]} & < & \infty\;, \quad p\geq 2\;,
\enq
and 
\beq\label{thomasse}
\sup_{\theta\in[0,T]}\Big(\| \alpha^1(\theta)\|_{H^p[\theta, T]} + \| \alpha^2(\theta)\|_{H^p[\theta, T]}\Big) &  < & \infty\;, \quad p\geq 2\;. 
\enq
We now write for $t \in [t_i, t_{i+1})$
\beqs
\E \big[ | Z_t(\theta) - Z_{t_i}(\theta) |^2 \big] & \leq & K ( I^1_{t_i, t} + I^2_{t_i, t} ) \;,
\enqs
with
\begin{equation*} \left\{
\begin{aligned}
  I^1_{t_i,t} & := ~ \E \big[ | \tilde Z_t(\theta) - \tilde Z _{t_i}(\theta) |^2 \big| \sigma\big(t_i, X^1_{t_i}(\theta)\big)\big| ^2 \big] \;,\\
  I^2_{t_i,t} & := ~ \E \big[ | \tilde Z_t(\theta)|^2 \big| \sigma\big(t, X_{t}^1(\theta)\big) - \sigma\big(t_i, X^1_{t_i}(\theta)\big)\big| ^2 \big] \;.
   \end{aligned}
\right.\end{equation*}
We give an upper bound for each term.
\beqs
  I^1_{t_i,t} & = & \E \Big[ \E \big[ | \tilde Z_t(\theta) - \tilde Z _{t_i}(\theta) |^2 \big| \Fc^1_{t_i}(\theta) \big]  \big| \sigma\big(t_i, X^1_{t_i}(\theta)\big) \big|^2 \Big] \\
 & \leq &   K\; \E \Big[ \int_{t_i}^{t_{i+1}} \big( | \alpha^1_r(\theta) |^2 + | \alpha^2_r(\theta) |^2 \big) dr \sup_{t\in[\theta,T]} \big| \sigma\big(t, X_{t}^1(\theta)\big) \big|^2  \Big] 
\enqs
which implies 
\beqs
\int_{t_i}^{t_{i+1}} I^1_{t_i,t} dt & \leq & K | \pi | \E \Big[ \int_{t_i}^{t_{i+1}}  \big( | \alpha^1_r (\theta)|^2 + | \alpha^2_r(\theta) |^2 \big) dr \sup_{t\in[\theta,T]} \big| \sigma\big(t, X_{t}^1(\theta)\big) \big|^2  \Big]  \;,
\enqs
therefore we have 
\beqs
\sum_{i=0,\; t_i \geq \theta}^{n-1} \int_{t_i}^{t_{i+1}} I^1_{t_i,t} dt & \leq &  K | \pi | \E \Big[ \int_{\theta}^{T}  \big( | \alpha^1_r(\theta) |^2 + | \alpha^2_r(\theta) |^2 \big) dr \sup_{t\in[\theta,T]} \big| \sigma\big(t, X_{t}^1(\theta)\big) \big|^2  \Big]  
 \;.
\enqs
From H\"older's inequality and \textbf{(HFD)} and \textbf{(HF)}, we have
\beqs
\sum_{i=0,\; t_i \geq \theta}^{n-1} \int_{t_i}^{t_{i+1}} I^1_{t_i,t} dt & \leq &  K | \pi | \E \Big[ \int_{\theta}^{T}  \big( | \alpha^1_r(\theta) |^4 + | \alpha^2_r |^4(\theta) \big) dr \Big]^{1\over 2} \Big(1+\E\Big[\sup_{t\in[\theta,T]}\big|  X_{t}^1(\theta) \big|^4 \Big] ^{1\over2}\Big)\;. 
\enqs
Using  \reff{thomasse}, we get
\beq\nonumber
\sum_{i=0, \;t_i \geq \theta}^{n-1} \int_{t_i}^{t_{i+1}} I^1_{t_i,t} dt & \leq &  K | \pi | \Big(1+\E\Big[\sup_{t\in[\theta,T]} \big|  X_{t}^1(\theta) \big|^4 \Big] ^{1\over2}\Big)  \\
 & \leq & K | \pi | \Big(1+ \E \Big[   \big|  X_{\theta}^1(\theta) \big|^4  \Big]^{1\over2}\Big) \label{inegI1theta}
 \;.
\enq 
We get from \reff{avantthomasse}, \textbf{(HFD)} and \textbf{(HF)}
\beqs
 I^2_{t_i,t} & \leq & 
  K \Big( \E \Big[ \big| \tilde Z_t - \tilde Z_{t_i} \big|^2 \big| X^1_{t_i} ( \theta ) \big|^2 \Big] + \E \Big[ \big| X^1_t(\theta) \tilde Z_t - X^1_{t_i}(\theta) \tilde Z_{t_i} \big|^2 \Big] +|\pi|^2\Big) \;.
\enqs
Arguing as above, we obtain
\beqs
\sum_{i=0,\; t_i \geq \theta}^{n-1} \int_{t_i}^{t_{i+1}} \E \Big[ \big| \tilde Z_t - \tilde Z_{t_i} \big|^2 \big| X^1_{t_i} ( \theta ) \big|^2 \Big]dt & \leq & K | \pi | \E \Big[ \sup_{\theta \leq t \leq T} \big( 1 + | X^1_t ( \theta ) |^4 \big) \Big]^{\frac{1}{2}} \;.
\enqs
Moreover, from It\^o's formula, $X^1(\theta) \tilde Z$ is a semimartingale of the form
\beqs
X^1_t(\theta) \tilde Z_t & = & X^1_\theta ( \theta ) \tilde Z_\theta + \int_\theta^t \tilde \alpha^1_r dr  + \int_\theta^t \tilde \alpha^2_r dW_r \;,
\enqs
where $|| \tilde \alpha^1||_{H^2[\theta, T]} + ||\tilde \alpha^2||_{H^2[\theta, T]}   \leq  K ( 1 + \E [ | X^1_{\theta} ( \theta) |^4 ] ^{1\over 4}\big)$. Therefore, we have
\beqs
\E \Big[ \big| X^1_t( \theta) \tilde Z_t - X^1_{t_i} ( \theta) \tilde{Z}_{t_i} \big|^2 \Big] & \leq &K \; \E \Big[ \int_{t_i}^{t_{i+1}} \big( | \tilde \alpha^1_r |^2 + | \tilde \alpha^2_r |^2 \big) dr \Big] \;,
\enqs
which implies
\beq
\sum_{i=0,~t_i\geq\theta}^{n-1} \int_{t_i}^{t_{i+1}} \E \Big[ \big| X^1_t(\theta) \tilde Z_t - X^1_{t_i}(\theta) \tilde Z_{t_i} \big|^2 \Big] & \leq &  K | \pi | \E \Big[ \big( 1 + | X^1_\theta ( \theta ) |^4 \big) \Big]^{\frac{1}{2}}   \;.\label{regXZtilde}
\enq
Using \reff{inegI1theta} and \reff{regXZtilde} we get the result.

\vspace{2mm}

When $b$, $\sigma$, $\beta$,  $f$ and $g$ are not in $C^1_b$, we can also prove the result by regularization. We first suppose that $f$ and $g$ are in $C^1_b$. We consider a density $q$ which is $C^\infty_b$ on $\R$ with a compact support, and we define an approximation $(b^\epsilon, \sigma^\epsilon, \beta^\epsilon)$ of $(b, \sigma,\beta)$ in $C^1_b$ by
\beqs
(b^\epsilon, \sigma^\epsilon, \beta^\epsilon)(t, x) & = & \frac{1}{\eps}\int_{\R} (b, \sigma,\beta) (t, x') q \Big( \frac{x-x'}{\epsilon} \Big) dx'\;,\quad (t,x)\in[0,T]\times\R\;.
\enqs
We then use the convergence of $(X^{1, \epsilon} (\theta), Y^{1, \epsilon} (\theta), Z^{1, \epsilon} (\theta))$ to $(X^{1} (\theta), Y^{1} (\theta), Z^{1} (\theta))$ and we get the result.
Next we assume that $f$ and $g$ are not $C^1_b$ and we consider for that $f^\epsilon$ and $g^\epsilon$ which are defined as previously and we get the result.
\ep

\vspace{2mm}

Using the link between $X^0$ and $X^1_\theta(\theta)$, we obtain that the bound \reff{borne-reg-Z1} is actually uniform in $\theta$.
 \begin{Corollary}\label{regularite z1 unif}
 Under \textbf{(HF)}, \textbf{(HFD)}, \textbf{(HBL)} and \textbf{(HBLD)}, there exists a constant $K$ such that 
\beq\label{borne-reg-Z1theta}
\E \Big[ \int_\theta^T \big| Z^1_t(\theta) - Z^1_{\pi(t)}(\theta) \big|^2 dt \Big] & \leq & K  |\pi| \;,
\enq
for all $\theta \in \pi$.
 \end{Corollary}
\ni \textbf{Proof.}
Since $X^0$ is a Brownian diffusion, we have for any $p \geq 2$, from \textbf{(HFD)} and \textbf{(HF)}, that 
\beqs
\E \Big[ \sup_{t \in [0,T]} | X^0_t|^p \Big] < \infty\;.
\enqs
 We notice that from the Lipschitz property of $\beta$ we have 
\beqs
\E \Big[\big| X^1_\theta ( \theta) \big|^4\Big] & = & \E \Big[\big|X^0_\theta+\beta\big(\theta, X^0_\theta\big) \big|^4\Big]\\
 & \leq & K\Big( 1+\E\Big[\sup_{t\in[0,T]}\big|X^0_t\big|^4\Big] \Big)~<~\infty\;.
\enqs
Combining this result with \reff{borne-reg-Z1}, we get \reff{borne-reg-Z1theta}
\ep

\vspace{2mm}

We now study the regularity of $Z^0$. 
\begin{Proposition}\label{regularite z0}
Under \textbf{(HF)}, \textbf{(HFD)}, \textbf{(HBL)} and \textbf{(HBLD)}, there exists a constant $K$ such that we have
\beqs
\E \Big[  \int_0^T \big| Z^0_t - Z^0_{\pi(t)} \big|^2 dt \Big] & \leq & K |\pi| \;.
\enqs
\end{Proposition}

\ni\textbf{Proof.}
The proof is similar to the previous one. The only difference is that the BSDE \reff{BSDE0} involves $Y^1$. We denote $\Theta^0_r = (r, X^0_r, Y^0_r, Z^0_r, Y^1_r(r) - Y^0_r)$. We first suppose that  $b$, $\sigma$, $\beta$, $f$ and $g$ are in $C^1_b$. We recall that 
\beqs
Y^0_t  & = & g(X^0_T) + \int_t^T f \big(\Theta^0_s \big) ds - \int_t^T Z^0_s dW_s \,.
\enqs
Therefore, for $0 \leq r \leq t \leq T$, we have
\beqs
D_r Y^0_t & = & \nabla g ( X^0_T) D_r X^0_T + \int_t^T \Big( \partial_x f \big( \Theta^0_s\big) D_r X^0_s + ( \partial_y - \partial_u) f\big( \Theta^0_s\big) D_r Y^0_s\\
& & + \partial_z f\big(\Theta^0_s\big) D_r Z^0_s +  \partial_u f\big(\Theta^0_s\big) D_r Y^1_s(s) \Big) dr - \int_t^T D_r Z^0_s dW_s\;,
\enqs
where $D X^0_r$, $D Y^0_r$, $D Z^0_r$ and $D Y^1_r(r)$ denote the Malliavin derivatives of  $X^0_r$,  $Y^0_r$,  $Z^0_r$ and  $Y^1_r(r)$ for $r\in[0,T]$. 
Using Malliavin calculus, we obtain  that a version of $Z^0$ is given by $(D_t Y^0_t)_{t \in [0, T]}$. By It\^{o}'s formula, we get
\beqs
\Lambda_t M_t Z_t & = & \E \Big[ M_T \Big( \Lambda_T \nabla g\big(X^0_T\big) D_t X^0_T + \int_t^T \big( \partial_x f\big( \Theta^0_r\big) D_t X^0_r + \partial_u f\big( \Theta^0_r\big) D_t Y^1_r(r) \big) \Lambda_r dr \Big)\Big| \Fc^0_t \Big] \;,
\enqs
where $\Lambda_t  :=  \exp ( \int_0^t ( \partial_y - \partial_u) f( \Theta^0_r) dr )$ and $M_t := 1 + \int_0^t M_r \partial_z f ( \Theta^0_r) dW_r$. Denote by  $\nabla X^0_t := \frac{\partial X^0_t }{\partial X^0_0}$ and $\nabla X^1_t(\theta) := \frac{\partial X^1_t(\theta)}{\partial X^1_0(\theta)}$ for $0\leq t\leq \theta\leq T$.
We then have  for $r\leq s\leq T$ 
\beqs
D_rX_s^1(s) =  (1 + \partial_x\beta(s,X^0_s)) D_r X^0_s= (1+ \partial_x\beta(s,X^0_s)) \nabla  X^0_s \sigma(r,X_r^0)[\nabla  X^0_r]^{-1} \;,
\enqs
thus we can see that $D_rX_s^1(s) = \nabla  X^1_s(s)\sigma(r,X_r^0)[\nabla  X^0_r]^{-1}$.
Therefore, we get  by writing the SDEs satisfied by $(D_r X^1_s(\theta))_{s\in[\theta,T]}$ for $r\leq \theta$, and $(\nabla X^1_s(\theta))_{s\in[\theta,T]}$
\beqs
D_r X^1_s(\theta) & = & \nabla X^1_s(\theta) \big[ \nabla X^0_r \big]^{-1} \sigma\big(r, X^0_r\big)\;,\quad r\leq \theta\leq s \;.
\enqs
Writing the BSDEs satisfied by $(D_r Y^1_s(\theta))_{s\in[\theta,T]}$ for $r\leq \theta$ and $(\nabla Y^1_s(\theta))_{s\in[\theta,T]}$, and using the previous equality, we get
\beqs
D_r Y^1_s(s) & = & \nabla Y^1_s(s) \big[ \nabla X^0_r \big]^{-1} \sigma\big(r, X^0_r\big)\;,\quad s\leq \theta \;.
\enqs
This implies
\beqs
\Lambda_t M_t Z_t & = & \E \Big[ M_T \Big(\Lambda_T \nabla g\big(X^0_T\big) \nabla X^0_T + \int_t^T F_r \Lambda_r dr \Big)\Big] \big[ \nabla X^0_t \big]^{-1} \sigma\big(t, X^0_t\big) \;,
\enqs
with $F_r := \partial_x f (\Theta^0_r) \nabla X^0_r + \partial_u f(\Theta^0_r) \nabla Y^1_r(r)$. We can write
\beqs
\Lambda_t M_t Z_t & = & \Big( \E [ G | \Fc^0_t ] - \int_0^t M_t F_r \Lambda_r dr \Big) \big[ \nabla X^0_t \big]^{-1} \sigma\big(t, X^0_t\big) \;,
\enqs
with $G := M_T ( \Lambda_T \nabla g (X^0_T) \nabla X^0_T + \int_0^T F_r \Lambda_r dr)$. Since $b$, $\sigma$, $f$ and $g$ have bounded derivatives, we have 
\beq\label{momG0}
\E\big[|G|^p\big] & < & \infty\;, \quad p\geq 2\;. 
\enq
Define $m_r := \E [ G | \Fc^0_r ]$ for $r\in[0,T]$. From \reff{momG0}  and  Doob's inequality, we have
\beq\label{momm0}
\| m \|_{S^p[0,T]} & < & \infty\;, \quad p\geq 2\;. 
\enq
Hence, there exists a process $\phi$ such that 
\beqs
m_r &= & \E[G  ] + \int_0^r \phi_u dW_u \;,\quad r\in[0, T]\;,
\enqs
and 
\beqs
\| \phi \|_{H^p[0,T]} & < & \infty\;, \quad p\geq 2\;. 
\enqs
We define $\tilde Z$ by 
\beqs
\tilde Z_t &: = & (\Lambda_t M_t)^{-1} \Big( m_t - M_t \int_0^t F_r \Lambda_r dr \Big) \big[ \nabla X^0_t \big]^{-1} \;,\quad t\in[0,T]\,.
\enqs
By It\^{o}'s formula, we can write
\beqs
\tilde Z_t & = & \tilde Z_0 + \int_0^t \alpha^1_r ds + \int_0^t \alpha^2_r dW_r \;,\quad t\in[0,T]\;.
\enqs
Using the fact that  $b$, $\sigma$, $f$ and $g$ have bounded derivatives and \reff{momm0}, we get 
\beqs
|| \tilde Z ||^p_{\Sc^p[0,T]} & < & \infty \;, \quad p \geq 2 \;,
\enqs
and 
\beq\label{idris}
\| \alpha^1\|_{H^p[0, T]} + \| \alpha^2\|_{H^p[0, T]} &  < & \infty\;, \quad p\geq 2\;. 
\enq
We now write for $t \in [t_i, t_{i+1})$
\beqs
\E \big[ | Z^0_t - Z^0_{t_i} |^2 \big] & \leq & K ( I^1_{t_i, t} + I^2_{t_i, t} ) \;,
\enqs
with 
\begin{equation*} \left\{
\begin{aligned}
  I^1_{t_i,t} & := ~ \E \big[ | \tilde Z_t - \tilde Z _{t_i} |^2 | \sigma(t_i, X^0_{t_i})| ^2 \big] \;,\\
  I^2_{t_i,t} & := ~ \E \big[ | \tilde Z_t|^2 \big| \sigma\big(t, X_{t}^0\big) - \sigma\big(t_i, X^0_{t_i}\big)\big| ^2 \big] \;.
   \end{aligned}
\right.\end{equation*}
As previously we give an upper bound for each term.
\beqs
I^1_{t_i, t} & \leq & K \; \E \Big[ \int_{t_i}^{t_{i+1}} \big( | \alpha^1_r |^2 + | \alpha^2_r |^2 \big) dr\sup_{t\in[0,T]} \big| \sigma\big(t, X_{t}^0\big) \big|^2  \Big] \;.
\enqs
From H\"older's inequality and Lipschitz property of $\sigma$, we have
\beqs
\sum_{i=0}^{n-1} \int_{t_i}^{t_{i+1}} I^1_{t_i,t} dt & \leq &  K | \pi | \E \Big[ \int_{0}^{T}  \big( | \alpha^1_r |^4 + | \alpha^2_r |^4 \big) dr \Big]^{1\over 2} \Big(1+\E\Big[\sup_{t\in[0,T]}\big|  X^0_{t} \big|^4 \Big] ^{1\over2}\Big)\;. 
\enqs
Using \reff{idris}, we get 
\beqs
\sum_{i=0}^{n-1} \int_{t_i}^{t_{i+1}} I^1_{t_i,t} dt & \leq & K | \pi | \;.
\enqs 
From \textbf{(HFD)} and \textbf{(HF)}, we get 
\beqs
 I^2_{t_i,t} & \leq & K \Big( \E \Big[ \big| \tilde Z_t - \tilde Z_{t_i} \big|^2 \big| X^0_{t_i}  \big|^2 \Big] + \E \Big[ \big| X^0_t \tilde Z_t - X^0_{t_i}
  \tilde Z_{t_i} \big|^2 \Big] +|\pi|^2\Big) \;.
\enqs
Arguing as above, we obtain
\beqs
\sum_{i=0}^{n-1} \int_{t_i}^{t_{i+1}} \E \Big[ \big| \tilde Z_t - \tilde Z_{t_i} \big|^2 \big| X^0_{t_i}  \big|^2 \Big]dt & \leq & K | \pi | \Big(1+\E \Big[ \sup_{t \in [0, T]} \big| X^0_t \big|^4  \Big]^{\frac{1}{2}}\Big) \;.
\enqs
Moreover, $X^0 \tilde Z$ is a semimartingale of the form
\beqs
X^0_t \tilde Z_t & = & X^0_0 \tilde Z_0 + \int_0^t \tilde \alpha^1_r dr  + \int_0^t \tilde \alpha^2_r dW_r
\enqs
where $|| \tilde \alpha^1||_{H^2[0, T]} + ||\tilde \alpha^2||_{H^2[0, T]}   \leq  K$ and we have
\beqs
\E \Big[ \big| X^0_t \tilde Z_t - X^0_{t_i} \tilde{Z}_{t_i} \big|^2 \Big] & \leq &K \, \E \int_{t_i}^{t_{i+1}} \big( | \tilde \alpha^1_r |^2 + | \tilde \alpha^2_r |^2 \big) dr \;,
\enqs
which implies
\beqs
\sum_{i=0}^{n-1} \int_{t_i}^{t_{i+1}} \E \Big[ \big| X^0_t \tilde Z_t - X^0_{t_i} \tilde Z_{t_i} \big|^2 \Big] & \leq &  K | \pi |    \;.
\enqs
When $b$, $\sigma$, $f$ and $g$ are not $C^1_b$, we can also prove the result by regularization as for Proposition \ref{regularite z1}.
\ep

\subsection{Error estimates for the recursive system of BSDEs}

\ni We first state an estimate of the approximation error for $(Y^1,Z^1)$.

\begin{Proposition}  \label{erreur Y1 Ypi}
Under \textbf{(HF)}, \textbf{(HFD)}, \textbf{(HBL)} and \textbf{(HBLD)}, we have the following  estimate  
\beqs
\sup_{\theta \in[ 0 , T]}\Big\{\sup_{t\in [\theta , T]} \E \Big[ \big| Y^1_{t}(\theta) - \YunPi_{\pi(t)}(\pi ( \theta )) \big|^2 \Big] + \E \Big[ \int_{\theta}^T   \big| Z^1_{s} ( \theta ) - \ZunPi_{\pi(s)} ( \pi ( \theta ) ) \big|^2 ds \Big] \Big\} & \leq  & K |\pi| \;, 
\enqs
for some constant $K$ which does not depend on $\pi$. 
\end{Proposition}
\ni\textbf{Proof.}
Fix $\theta\in[0,T]$ and $t\in[\theta,T]$. We then have 
\beq\nonumber
\E \Big[ \big| Y^1_{t}(\theta) - \YunPi_{\pi(t)}(\pi ( \theta )) \big|^2 \Big]  & \leq &  2\; \E \Big[ \big| Y^1_{t}(\theta ) - Y^1_{t}(\pi ( \theta )) \big|^2 \Big]\\
 & & \label{decomperrschemeY1} + \;2 \;\E \Big[ \big| Y^1_{t}(\pi ( \theta )) - \YunPi_{\pi(t)}(\pi ( \theta )) \big|^2 \Big] \;.
\enq
We study separately the two terms of right hand side.

Define $\delta X^1_t(\theta) := X^1_{t}(\theta) - X^1_{t} ( \pi ( \theta ))$, $\delta Y^1_t(\theta) := Y^1_{t}(\theta) - Y^1_{t} ( \pi ( \theta ))$  and $\delta Z^1_t(\theta) := Z^1_{t}(\theta) - Z^1_{t} ( \pi ( \theta ))$. 
Applying It\^{o}'s formula, we get
\beqs
|\delta Y^1_T(\theta)|^2 - |\delta Y^1_t(\theta)|^2 & = & 
2 \int_{t}^T\hspace{-2mm} \delta Y^1_s(\theta) \Big[ f\big(\Theta^1_s( \pi(\theta)) \big) - f\big(\Theta^1_s( \theta)\big) \Big]ds \\
&& + \;2 \int_{t}^T \delta Y^1_s(\theta) \delta Z^1_{s} ( \theta )  dW_{s} +  \int_{t}^T |\delta Z^1_s(\theta)|^2ds\;, 
\enqs
where $\Theta^1_s (\theta) := (s, X^1_{s} ( \theta ), Y^1_{s} ( \theta ), Z^1_{s} ( \theta ), 0)$.
From \textbf{(HBL)} and \textbf{(HBLD)}, we get  
\beqs
\E \big[ |\delta Y^1_t(\theta)|^2\big] & \leq &  K \Big( \E \big [ | \delta X^1_T (\theta)|^2 \big] +   \E \Big[ \int_{t}^T | \delta Y^1_s(\theta)| |\delta X^1_s (\theta)| ds \Big] +  \E \Big[  \int_t^T  \hspace{-1mm}| \delta Y^1_s(\theta)|^2 ds \Big]\\
& &  
 +   \; \E \Big[ \int_{t}^T \hspace{-1mm}  | \delta Y^1_{s}(\theta) | |\delta Z^1_{s}(\theta) | ds \Big] \Big)
  \hspace{-1mm}  -  \E \Big[ \int_{t}^T \hspace{-1mm}  |\delta Z^1_{s} ( \theta )|^2 ds \Big]\;. 
\enqs
Using the inequality $2 ab \leq a^2 / \eta + \eta b^2$ for $a,b\in\R$ and $\eta>0$, we can see that
\beq\nonumber
\E \big[ |\delta Y^1_t(\theta)|^2\big ]  +  \E \Big[ \int_{t}^T |\delta Z^1_{s} ( \theta )|^2 ds \Big]& \leq & K \Big( \E \big[ | \delta X^1_T (\theta)|^2 \big]   +  \int_{t}^T \E \big[ |\delta Y^1_s(\theta)|^2\big ]  ds \\
& &  + \; \E \Big[ \int_{t}^T | \delta X^1_s (\theta)|^2 ds \Big]  \Big)\;. \label{majY1Z1Gronwall}
\enq 
From \reff{maj1DX1} and Gronwall's lemma, we get 
\beq\label{majdeltaY1(theta)}
\E \big[ \big|Y^1_{t}(\theta) - Y^1_{t} ( \pi ( \theta )) \big|^2 \big]  & \leq & K |\pi| \;.
\enq

\vspace{2mm}

We now study the  second term of the right hand side of \reff{decomperrschemeY1}. Using the same argument as in the proof of Theorem 3.1 in \cite{boutou04}, we get from the regularity of $Z^1$ given by Corollary \ref{regularite z1 unif}  
\beq\label{etoile}
 \E \Big[ \big| Y^1_{t}(\pi ( \theta )) - \YunPi_{\pi(t)}(\pi ( \theta )) \big|^2 \Big]  & \leq &  K|\pi|\;.
\enq
This last inequality with \reff{decomperrschemeY1} and  \reff{majdeltaY1(theta)} gives
\beqs
\sup_{\theta\in[0,T]} \Big\{\sup_{t\in[\theta,T]}\E\Big[ \big| Y^1_{t}(\theta) - \YunPi_{\pi(t)}(\pi ( \theta )) \big|^2 \Big] \Big\} & \leq & K |\pi|\;.
\enqs
\ni We now turn to the error on the term $Z^1(\theta)$. We first use the inequality
\beq\nonumber
\E \Big[  \int_{\theta}^T \big| Z^1_{t}(\theta) - \ZunPi_{\pi(t)} ( \pi(\theta)) \big|^2 dt \Big]  & \leq &  
2\; \E \Big[  \int_{\theta}^T \big| Z^1_{t} ( \pi(\theta)) - \ZunPi_{\pi(t)} ( \pi(\theta)) \big|^2 dt \Big]\\
& &\label{decomperrZ1}
 + \;2 \;\E \Big[  \int_{\theta}^T \big| \delta Z^1_{t}(\theta)  \big|^2 dt \Big]  \;.
\enq
Using \reff{majY1Z1Gronwall} and \reff{majdeltaY1(theta)} with $t=\theta$, we get
\beq\label{errdeltaZ1}
\E \Big[  \int_{\theta}^T \big|\delta Z^1_{s}(\theta)  \big|^2 ds \Big] & \leq & K |\pi| \;.
\enq
The other term in the right hand side of \reff{decomperrZ1} is the classical error in an approximation of BSDE. Therefore, using Corollary \ref{regularite z1 unif}  and \reff{etoile}, we have 
\beq\label{errdeltaZ1pitheta}
\E \Big[  \int_{\theta}^T \big| Z^1_{t} ( \pi(\theta)) - \ZunPi_{\pi(t)} ( \pi(\theta)) \big|^2 dt \Big] & \leq & K |\pi| \;.
\enq
Combining \reff{decomperrZ1}, \reff{errdeltaZ1} and  \reff{errdeltaZ1pitheta}, we get
\beqs
\E \Big[  \int_{\theta}^T \big| Z^1_{t} ( \theta) - \ZunPi_{\pi(t)} ( \pi(\theta)) \big|^2 dt \Big] & \leq & K |\pi|\;.
\enqs
\ep

\vspace{2mm}

We now turn to  the estimation of the error between $(Y^0,Z^0)$ and its Euler scheme \reff{YZ0}. 
Since this scheme involves the approximation $Y^{1,\pi}$ of $Y^1$, we first need to introduce an intermediary scheme involving the "true" value of the process $Y^1$. 
We therefore consider the scheme $(\YtilzePi,\ZtilzePi)$ defined by 
\begin{equation}\label{YZ0bis} \left\{
\begin{aligned}
  \YtilzePi_{T} ~ & = ~ g ( \XzePi_{T} ) \;,\\
  \YtilzePi_{t_{i-1}} ~ & = ~ \E^0_{i-1} \big[ \YtilzePi_{t_{i}} \big] + f \big( t_{i-1},\XzePi_{t_{i-1}}, \YtilzePi_{t_{i-1}}, \ZtilzePi_{t_{i-1}}, Y^1_{t_{i-1}}(t_{i-1}) - \YtilzePi_{t_{i-1}} \big) \Delta t^\pi_{i} \;,\\
    \ZtilzePi_{t_{i-1}} ~ & = ~ \frac{1}{\Delta t^\pi_{i} } \E^0_{i-1} \big[ \YtilzePi_{t_i} \Delta W^\pi_{i} \big] \;,\quad1\leq i\leq n\;.
    \end{aligned}
\right.\end{equation}
Using the regularity result of Proposition \ref{regularite z0} and the same arguments as in the proof of Theorem 3.1 in \cite{boutou04},  we get under \textbf{(HF)}, \textbf{(HFD)}, \textbf{(HBL)} and \textbf{(HBLD)}
\beq\label{majerrtilde0}
\sup_{t\in[0,T]} \E \Big[ \big| Y^0_{t} - \YtilzePi_{\pi ( t )}\big|^2 \Big] + \E \Big[ \int_{0}^T \big|Z^0_{t} - \ZtilzePi_{\pi (t)}\big|^2 dt \Big] & \leq & K |\pi|\;.
\enq
With this inequality, we get the following estimate for the error between $(Y^0, Z^0)$ and the Euler scheme \reff{YZ0}.
\begin{Proposition}\label{erreur Y0 Ypi}
Under \textbf{(HF)}, \textbf{(HFD)}, \textbf{(HBL)} and \textbf{(HBLD)},
we have the following estimate
\beqs
\sup_{t\in[0,T]} \E \Big[ \big| Y^0_{t} - \YzePi_{\pi ( t )}\big|^2 \Big] + \E \Big[ \int_{0}^T \big|Z^0_{t} - \ZzePi_{\pi (t)}\big|^2 dt \Big] & \leq & K |\pi|\;,
\enqs
for some constant $K$ which does not depend on $\pi$.
\end{Proposition}

\ni\textbf{Proof.} 
We first remark that
\begin{equation*}
\left\{ \begin{aligned}
\sup_{t\in[0,T]} \E \Big[ \big| Y^0_{t} - \YzePi_{\pi ( t )}\big|^2 \Big] & \leq & 2 \sup_{t\in[0,T]} \E \Big[ \big| Y^0_{t} - \YtilzePi_{\pi ( t )}\big|^2 \Big] + 2 \sup_{t\in[0,T]} \E \Big[ \big| \YzePi_{\pi ( t )} - \YtilzePi_{\pi ( t )}\big|^2 \Big] \;,\\
 \E \Big[ \int_{0}^T \big|Z^0_{t} - \ZzePi_{\pi (t)}\big|^2 dt \Big] & \leq &  2 \;\E \Big[ \int_{0}^T \big|Z^0_{t} - \ZtilzePi_{\pi (t)}\big|^2 dt \Big] +  2 \;\E \Big[ \int_{0}^T \big|\ZzePi_{\pi (t)} - \ZtilzePi_{\pi (t)} \big|^2 dt \Big] \;.
\end{aligned}\right.
\end{equation*}
Using \reff{majerrtilde0}, we only need to study $\sup_{t\in[0,T]} \E [ | \YzePi_{\pi ( t )} - \YtilzePi_{\pi ( t )}|^2 ]$ and $\E [ \int_{0}^T |\ZzePi_{\pi (t)} - \ZtilzePi_{\pi (t)} |^2 dt ]$. To this end, we need to introduce continuous schemes for all $0 \leq i \leq n-1$. Since $\E[|\YzePi_{t_i}|^2] < \infty$ and $\E[|\YtilzePi_{t_i}|^2] < \infty$  for all $1 \leq i \leq n$, we deduce, from the martingale representation theorem, that there exist square integrable processes $\uZzePi$ and $ \uZtilzePi$ such that
\beqs
\YzePi_{t_i} &=& \E \big[ \YzePi_{t_{i+1}} \big| \Fc_{t_i} \big] + \int_{t_i}^{t_{i+1}}  \uZzePi_s dW_s \;,\\
\YtilzePi_{t_i} &=& \E \big[ \YtilzePi_{t_{i+1}} \big| \Fc_{t_i} \big] + \int_{t_i}^{t_{i+1}}  \uZtilzePi_s dW_s \;.
\enqs
We then define 
\begin{equation*}
\left\{ \begin{aligned}
\YzePi_t & = ~ \YzePi_{t_i} - (t - t_i) f\big(t_i, \XzePi_{t_i}, \YzePi_{t_i}, \ZzePi_{t_i}, \YunPi_{t_i}(t_i) - \YzePi_{t_i}\big) + \int_{t_i}^t \uZzePi_s dW_s \;,\\
\YtilzePi_t & = ~ \YtilzePi_{t_i} - (t - t_i) f\big(t_i, \XzePi_{t_i}, \YtilzePi_{t_i}, \ZtilzePi_{t_i}, Y^1_{t_i}(t_i) - \YtilzePi_{t_i}\big) + \int_{t_i}^t  \uZtilzePi_s dW_s \;,
\end{aligned}\right.
\end{equation*}
for $t\in[t_i,t_{i+1})$. 
Let $i \in \{0, \ldots, n-1\}$ be fixed, and set $\delta Y_t := \YzePi_t - \YtilzePi_t $, $\delta Z_i:= \ZzePi_{t_i} - \ZtilzePi_{t_i}$, $\delta \underline Z_t := \uZzePi_t - \uZtilzePi_t$ and $\delta f_t : = f(t_i, \XzePi_{t_i}, \YzePi_{t_i}, \ZzePi_{t_i}, \YunPi_{t_i}(t_i) - \YzePi_{t_i}) -  f(t_i, \XzePi_{t_i}, \YtilzePi_{t_i}, \ZtilzePi_{t_i}, Y^1_{t_i}(t_i) - \YtilzePi_{t_i})$ for $t \in [t_i, t_{i+1})$. By It\^{o}'s formula, we compute that
\beqs
A_t & := & \E | \delta Y_t |^2 + \int_t^{t_{i+1}} \E | \delta \underline Z_s |^2 ds  - \E | \delta Y_{t_{i+1}} |^2 ~ = ~2 \int_t^{t_{i+1}} \E [ \delta Y_s \delta f_s ] ds \;, \quad t_i \leq t \leq t_{i+1} \;.
\enqs
Let $\alpha > 0$ be a constant to be chosen later on. From the Lipschitz property of $f$ and the inequality $2ab \leq \alpha a^2 + b^2 / \alpha$, we get
\beqs
A_t & \leq & \alpha \int_t^{t_{i+1}} \E | \delta Y_s |^2 ds + \frac{K}{\alpha} \int_t^{t_{i+1}} \E \Big[  | \delta Y_{t_i} |^2 + | \delta Z_{i} |^2 + | Y^1_{t_i}(t_i) - \YunPi_{t_i}(t_i) |^2 \Big] ds \;.
\enqs
Using Proposition \ref{erreur Y1 Ypi}, we get
\beqs
A_t & \leq & \alpha \int_t^{t_{i+1}} \E | \delta Y_s |^2 ds + \frac{K}{\alpha} |\pi| \; \E | \delta Y_{t_i} |^2 + \frac{K}{\alpha} \int_t^{t_{i+1}}\E | \delta Z_{i} |^2 ds + \frac{K}{\alpha} |\pi|^2 \;.
\enqs
We can write
\beq \label{double inegalite}
\E | \delta Y_t|^2 & \leq & \E | \delta Y_{t_{i+1}}|^2 + \int_t^{t_{i+1}} \E | \delta \underline Z_s |^2 ds ~ \leq ~ \alpha \int_t^{t_{i+1}} \E | \delta Y_s |^2 ds + B_i \:,
\enq
where
\beqs
B_i & := & \E | \delta Y_{t_{i+1}} |^2 + \frac{K}{\alpha} | \pi|\;  \E | \delta Z_{i} |^2  + \frac{K}{\alpha} |\pi| \; \E | \delta Y_{t_i} |^2  + \frac{K}{\alpha} |\pi|^2 \;.
\enqs
By Gronwall's lemma, this shows that $\E | \delta Y_t |^2 \leq B_i e^{\alpha | \pi |}$ for $t_i \leq t < t_{i+1}$, which plugged in the second inequality of \reff{double inegalite} provides
\beq\label{inegalite gronwall discret}
 \E | \delta Y_t|^2 +  \int_t^{t_{i+1}} \E | \delta \underline Z_s |^2 ds  & \leq & B_i \Big( 1 + \alpha | \pi | e^{\alpha |\pi|} \Big)\;. 
\enq
Interpreting $Z^{0,\pi}_{t_i}$ (resp.  $ \tilde Z_{t_i}^{0,\pi}$) as the projection  of $\underline Z^{0,\pi}$ (resp.  $\underline {\tilde Z}^{0,\pi}$)  in $H^2_\F[t_i,t_{i+1}]$ on the set of constant processes, we have
\beq\label{inegprojZ}
\int_{t_i}^{t_{i+1}}\E | \delta Z_{i} |^2 ds &  \leq & \int_{t_i}^{t_{i+1}}\E | \delta \underline Z_s |^2 ds\;.
\enq 
Applying \reff{inegalite gronwall discret} for $t = t_i$ and $\alpha=2K$, and using the previous inequality, we get
\beqs
 \E | \delta Y_{t_i} |^2 +k_1(\pi) \int_{t_i}^{t_{i+1}} \E | \delta \underline Z_s |^2 ds  & \leq & k_2(\pi) \E | \delta Y_{t_{i+1}} |^2 + k_3(\pi)
| \pi |^2 \;,\quad 0\leq i\leq n-1\;,
\enqs
where $k_1(\pi) = \frac{{1\over 2}-  K|\pi|e^{2K|\pi|}}{1 -  {|\pi|\over 2} -  K | \pi|^2e^{2K|\pi|}}$, $k_2(\pi) = \frac{1 + 2 K|\pi|e^{2K|\pi|}}{1 - { |\pi|\over 2} -  K | \pi|^2e^{2K|\pi|}}$ and $k_3(\pi) = \frac{{1\over 2} +  K|\pi|e^{2K|\pi|}}{1 - { |\pi|\over 2} -  K | \pi|^2e^{2K|\pi|}}$.
Since for small $|\pi|$ we have $k_1(\pi)$ $\geq$ $0$, we get
\beqs
 \E | \delta Y_{t_i} |^2  & \leq & k_2(\pi) \E | \delta Y_{t_{i+1}} |^2 + k_3(\pi)
| \pi |^2 \;,\quad 0\leq i\leq n-1\;,
\enqs
 for $|\pi|$ small enough. 
 
\ni Iterating this inequality, we get
\beqs
\E | \delta Y_{t_i} |^2   & \leq & k_2(\pi)^{1\over |\pi|} \E | \delta Y_{t_n} |^2 + | \pi |^2 k_3(\pi) \sum_{j={i}}^n k_2(\pi) ^{j-i}\;.
\enqs
Since $k_2 (\pi) \geq 1$ and $\delta Y_{t_n} = 0$, we get for small $|\pi|$
\beq\label{ecart sur i}
\E | \delta Y_{t_i} |^2 & \leq & | \pi | k_3(\pi) k_2(\pi) ^{1 \over |\pi|} ~\leq ~ K |\pi|\;,\quad 0\leq i\leq n \;,
\enq
which gives 
\beqs
\sup_{t\in[0,T]} \E \Big[ \big| \YzePi_{\pi ( t )} - \YtilzePi_{\pi ( t )}\big|^2 \Big] & \leq &  K|\pi|\;.
\enqs
 Summing up the inequality \reff{inegalite gronwall discret} with $t = t_i$ and $\alpha=2K$ and using \reff{inegprojZ}, we get 
 \beqs
\Big({1\over 2}-K|\pi|e^{2K|\pi|}\Big)\int_0^T \E |  Z^{0,\pi}_{\pi(s)}- \tilde Z^{0,\pi}_{\pi(s)}|^2 ds & \leq & 
2K|\pi|e^{2K|\pi|}\sum_{i=1}^{n-1}\E | \delta Y_{t_{i}} |^2+(1+2K|\pi|)\E | \delta Y_{t_{n}} |^2\\
 & & +\Big(\frac{1}{2}+K|\pi|e^{2K|\pi|}\Big)\Big(|\pi| + |\pi|\sum_{i=0}^{n-1}\E | \delta Y_{t_{i}} |^2 \Big)\;.
 \enqs
 Using \reff{ecart sur i}, we get for $|\pi|$ small enough
\beqs \label{ecart sur z}
\int_0^T \E |  Z^{0,\pi}_{\pi(s)}- \tilde Z^{0,\pi}_{\pi(s)}|^2 ds  & \leq & K | \pi|  \;.
\enqs
\ep

\subsection{Error estimate for the BSDE with a jump} 
We now give an error estimate of the approximation scheme for the BSDE with a jump.
\begin{Theorem} \label{erreur schema} Under \textbf{(HF)}, \textbf{(HFD)}, \textbf{(HBL)} and \textbf{(HBLD)},
we have the following  error estimate for the approximation scheme  
\beqs
\sup_{t\in[0,T] } \E \Big[ \big| Y_{t} - Y^{\pi}_{t} \big|^2 \Big] + \E \Big[ \int_0^T \big| Z_t - Z^\pi_{t} \big|^2 dt \Big] + \E \Big[ \int_0^T  \lambda_t  \big| U_t - U^\pi_{t} \big|^2 dt \Big] & \leq & K |\pi| \;,
\enqs
for some  constant $K$ which does not depend on $\pi$.
\end{Theorem}

\ni\textbf{Proof.} 

\ni\textbf{Step 1.} \textit{Error for the variable $Y$}. Fix $t\in[0,T]$. From Theorem \ref{theoreme existence unicite} and \reff{SchemeYZ}, we have 
\beqs
\E \Big[ \big| Y_{t} - Y^{\pi}_{t} \big|^2 \Big]  & = & \E \Big[ \big| Y^0_{t} - \YzePi_{\pi(t)} \big|^2\1_{t<\tau} \Big] +\E \Big[ \big| Y^1_{t}(\tau) - \YunPi_{\pi(t)}(\pi(\tau)) \big|^2\1_{t\geq\tau} \Big] \;.
\enqs
Using \textbf{(DH)}, we get
\beqs
\E \Big[ \big| Y_{t} - Y^{\pi}_{t} \big|^2 \Big]  & \leq & \E \Big[ \big| Y^0_{t} - \YzePi_{\pi(t)} \big|^2\Big] +\int_0^T\E \Big[ \big| Y^1_{t}(\theta) - \YunPi_{\pi(t)}(\pi(\theta)) \big|^2\1_{t\geq\theta} \gamma_T(\theta)\Big]d\theta \\
 & \leq & K\Big(\E \Big[ \big| Y^0_{t} - \YzePi_{\pi(t)} \big|^2\Big] + \sup_{\theta\in[0,T]} \sup_{s\in[\theta,T]}\E \Big[ \big| Y^1_{s}(\theta) - \YunPi_{\pi(s)}(\pi(\theta)) \big|^2\Big]\Big)\;.
 \enqs
 Using Propositions \ref{erreur Y1 Ypi} and \ref{erreur Y0 Ypi}, and  since $t$ is arbitrary chosen in $[0,T]$, we get
 \beqs
\sup_{t\in[0,T]}\E \Big[ \big| Y_{t} - Y^{\pi}_{t} \big|^2 \Big]   & \leq & K |\pi|\;.
\enqs

\ni\textbf{Step 2.} \textit{Error estimate for the variable $Z$}.
From Theorem \ref{theoreme existence unicite} and \reff{SchemeYZ}, we have 
\beqs
\E\Big[\int_0^{T} {\big|Z_t-Z^\pi_t\big|}^2dt\Big]  =  \E\Big[\int_0^{T\wedge\tau}{\big|Z^0_t-\ZzePi_{\pi(t)}\big|}^2dt\Big] +  \E\Big[\int_{T \wedge \tau}^T{\big|Z^1_t(\tau)-\ZunPi_{\pi(t)}(\pi(\tau))\big|}^2dt\Big]\;.
\enqs
Using \textbf{(DH)}, we get 
\beqs
\E\Big[\int_0^T{\big|Z_t-Z^\pi_t\big|}^2dt\Big] & = &  \hspace{-2mm} \int_0^T\int_0^\theta\E\Big[{\big|Z^0_t-\ZzePi_{\pi(t)}\big|}^2\gamma_T(\theta)\Big]dtd\theta \\
 & &  \hspace{-2mm} +  \int_0^T\int_\theta^T\E\Big[{\big|Z^1_t(\theta)-\ZunPi_{\pi(t)}(\pi(\theta))\big|}^2\gamma_T(\theta)\Big]dtd\theta\;. \\
 & \leq & \hspace{-2mm}  K\Big( \E\Big[\int_0^T \hspace{-1mm} {\big|Z^0_t-\ZzePi_{\pi(t)}\big|}^2dt\Big] +\sup_{\theta\in[0,T]}  \E\Big[\int_\theta^T \hspace{-1mm}{\big|Z^1_t(\theta)-\ZunPi_{\pi(t)}(\pi(\theta))\big|}^2\Big]dt\Big)\;.
\enqs
From Propositions \ref{erreur Y1 Ypi} and \ref{erreur Y0 Ypi}, we get
\beqs
\E\Big[\int_0^T{\big|Z_t-Z^\pi_t\big|}^2dt\Big] & \leq & K|\pi|\;.
\enqs

\ni\textbf{Step 3.} \textit{Error estimate for the variable $U$}.
From Theorem \ref{theoreme existence unicite} and \reff{SchemeYZ}, we have
\beqs
\E\Big[\int_0^T \big| U_t - U^\pi_{t} \big|^2\lambda_tdt\Big]  &\leq & K \, \E\Big[\int_0^T\Big(| Y^1_t(t) - \YunPi_{\pi(t)} ( \pi(t) ) |^2 + | Y^0_t - \YzePi_{\pi(t)}  |^2 \Big)\lambda_t dt\Big]\;.
\enqs
Using \textbf{(HBI)}, we get
\beqs
\E\Big[ \hspace{-1mm} \int_0^T \hspace{-1mm}  \big| U_t - U^\pi_{t} \big|^2\lambda_tdt\Big]   \leq  K\Big(\sup_{\theta\in[0,T]}\sup_{t\in[\theta,T]}\E\Big[\big| Y^1_t(\theta) - \YunPi_{\pi(t)} ( \pi(\theta) ) \big|^2\Big]+\sup_{t\in[0,T]}\E\Big[\big| Y^0_t - \YzePi_{\pi(t)}  \big|^2\Big] \Big)\,.
\enqs
Combining this last inequality with Propositions \ref{erreur Y1 Ypi} and \ref{erreur Y0 Ypi}, we get the result.
\ep
\begin{Remark}
{\rm Our decomposition approach allows us to suppose that the jump coefficient $\beta$ is only Lipschitz continuous. We do not need to impose any regularity condition on $\beta$ and any ellipticity assumption on $I_d+\nabla \beta$ as done in \cite{boueli08} in the case of Poissonian jumps independent of $W$.
}
\end{Remark}

\section{Convergence of the backward scheme for the quadratic case}
\setcounter{equation}{0} \setcounter{Assumption}{0}
\setcounter{Theorem}{0} \setcounter{Proposition}{0}
\setcounter{Corollary}{0} \setcounter{Lemma}{0}
\setcounter{Definition}{0} \setcounter{Remark}{0}

In this section we assume that \textbf{(HBQ)} holds and that $\sigma(t,x)=\sigma(t,0)=\sigma(t)$ for any $t \in \R_+$ and $x\in \R$.

Before giving the error of the scheme we give a uniform bound for the processes $Z^0$ and $Z^1$ which allows to prove that the BSDE \reff{BSDE} is Lipschitz and thus we can use Theorem \ref{erreur schema}. For that we introduce the BMO-martingales class, and we also give some bounds for the processes $X^0$, $X^1$, $Y^0$ and $Y^1$.

\subsection{BMO property for the solution of the BSDE}

To obtain a uniform bound for the processes $Z^0$ and $Z^1$ we need the following assumption.

%

\vspace{2mm}

\ni\textbf{(HBQD)} There exists a constant $K_f$ such that the function $f$ satisfies
\beqs
| f(t,  x, y, z, u) - f(t',  x', y', z', u') | & \leq & K_{f} \big[ |x - x'| + |y - y'| + |u - u'|+ |t - t'|^{1\over 2} \big]  \\ & &+ L_{f, z} ( 1+|z| + |z'| ) |z - z'| \;,
\enqs
for all $(t,x,x',y,y',z,z',u,u')\in [0,T] \times \R^2\times \R^2\times \R^2\times \R^2$.\\

\vspace{2mm}

In the sequel of this section, the space of BMO martingales plays a key role for the a priori estimates of processes $Z^0$ and $Z^1$. We refer to \cite{kaz94} for the theory of BMO martingales. Here, we just give the definition of a BMO martingale and recall a property that we use in the sequel. 
\begin{Definition}
A process $M$ is said to be a $BMO_\F[0,T]$-martingale if 
 $M$ is a square integrable $\F$-martingale s.t. 
\beqs
{\|M\|}_{BMO_\F[0,T]} & := & \sup_{\tau \in \Tc_\F[0,T]} \E \Big[ \big|M_T-M_\tau \big|^2 \Big|\Fc_\tau \Big]^{1/2} ~ < ~ \infty \;,
\enqs
where $\Tc_\F[0,T]$ denotes the set of $\F$-stopping times valued in $[0,T]$.
\end{Definition}
\ni The BMO condition provides a property on the Dolean-Dade exponential of the process $M$. 
\begin{Lemma}\label{ppte-gen-BMO}
Let $M$  be a $BMO_\F[0,T]$-martingale. Then the stochastic exponential $\Ec(M)$ defined by 
\beqs
\Ec ( M)_t & = & \exp \Big( M_t - \frac{1}{2} \langle M,M \rangle\Big)_t \;, \quad 0 \leq t \leq T\;,
\enqs
is a uniformly integrable $\F$-martingale.
\end{Lemma}
\noindent We refer to \cite{kaz94} for the proof of this result.\\

We first state a BMO property for the processes $Z^0$ and $Z^1$, which will be used in the sequel to provide an estimate for these processes.
\begin{Lemma}\label{BMO}
Under \textbf{(HF)}, \textbf{(HBQ)} and \textbf{(HBQD)}, 
the martingales $\int_0^.Z_s^0 d W_s$ and  $\int_0^.Z^1_s(\theta) \mathds{1}_{s\geq\theta}d W_s$, $\theta\in[0,T]$ are $BMO_\F[0,T]$-martingales and there exists a constant $K$ which is independent from $\theta$ such that
\beqs
{\Big\|\int_0^.Z^0_s d W_s\Big\|}_{BMO_\F[0,T]} & \leq & K\;,\\
\sup_{\theta\in[0,T]}{\Big\|\int_0^.Z^1_s(\theta) \mathds{1}_{s\geq\theta}d W_s\Big\|}_{BMO_\F[0,T]} & \leq & K\;.
\enqs
 \end{Lemma}

\ni\textbf{Proof.} 
Define the function $\phi:~\R\rightarrow\R$ by
\beq\label{ppte phi}
\phi(x) &  = & (e^{2K_qx} - 2K_qx - 1)/2 K_q^2\;,\quad x\in\R\;.
\enq
 We notice that $\phi$ satisfies 
 \beqs
 \phi'(x)  & \geq &  0 ~\mbox{ and }~ \frac{1}{2} \phi''(x) - K_q \phi'(x) ~=~ 1\;,
 \enqs
 for $x \geq 0$. Since $Y^0$ and $Y^1(.)$ are solutions to quadratic BSDEs with bounded terminal conditions, we get from Proposition 2.1 in \cite{k00}  
 the existence of a constant $m$ such that 
 \beq\label{born-unif-Y}
\| Y^0\|_{\Sc^\infty[0,T]}  & \leq &  m \quad \mbox{ and }\quad  \sup_{\theta\in[0,T]}\|Y^1(\theta)\|_{\Sc^\infty[\theta,T]}  ~ \leq ~  m\;.
\enq
Applying It\^o's formula we get 
\begin{multline*}
\phi(Y^0_\nu + m) + \E \Big( \int_\nu^T \frac{1}{2} \phi''(Y^0_s + m) | Z^0_s|^2 ds\Big| \Fc_\nu \Big)  =  \\
\E(\phi( Y^0_T + m) | \Fc_\nu) + \E \Big( \int_\nu^T \phi'(Y^0_s + m) f(s, X^0_s, Y^0_s, Z^0_s, Y^1_s(s) - Y^0_s) ds \Big| \Fc_\nu \Big) \;,
\end{multline*}
for any $\F$-stopping time $\nu$ valued in $[0,T]$. 
From the growth assumption on the generator $f$ in \textbf{(HBD)}, \reff{ppte phi} and \reff{born-unif-Y}, we obtain
\begin{multline*}
\phi(Y^0_\nu + m) + \E \Big( \int_\nu^T| Z^0_s|^2 ds\Big| \Fc_\nu \Big)  \leq \\
 \E(\phi( Y^0_T + m) | \Fc_\nu)
 + \E \Big( \int_\nu^T \phi'(Y^0_s + m) K_q (1 + 2 {\|Y^0\|}_{\Sc^{\infty}} + \sup_{\theta\in[0,T]}\|Y^1 (\theta)\|_{\Sc^{\infty}[\theta,T]}) ds \Big| \Fc_\nu \Big) \;. 
\end{multline*}
This last inequality and \reff{born-unif-Y} imply that there exists a constant $K$ which depends only on $m$, $T$ and $K_q$ such that for all $\F$-stopping times $\nu \in[0, T]$
\beqs
 \E \Big( \int_\nu^T| Z^0_s|^2 ds\Big| \Fc_\nu \Big) & \leq & K \;.
\enqs
For the process $Z^1$, we use the same technics. Let us fix $\theta \in [0, T]$. 
Applying It\^o's fomula we get
\begin{multline*}
\phi(Y^1_{\nu\vee \theta}(\theta) + m) + \E \Big( \int_{\nu\vee \theta}^T \frac{1}{2} \phi''(Y^1_s(\theta) + m) | Z^1_s(\theta)|^2ds \Big| \Fc_{\nu \vee \theta} \Big)  =   \\
 \E(\phi( Y^1_T(\theta) + m) | \Fc_\nu) + \E \Big( \int_{\nu \vee \theta}^T \phi'(Y^1_s(\theta) + m) f(s, X^1_s(\theta), Y^1_s(\theta), Z^1_s(\theta), 0) ds \Big| \Fc_{\nu \vee \theta} \Big) \;,
\end{multline*}
for any $\F$-stopping time $\nu$ valued in $[0,T]$.
From the growth assumption on the generator $f$ in \textbf{(HBQ)}, \reff{ppte phi} and \reff{born-unif-Y}, we obtain
\begin{multline*}
\phi(Y^1_{\nu \vee \theta}(\theta) + m) + \E \Big( \int_{\nu \vee \theta}^T| Z^1_s(\theta)|^2ds \Big| \Fc_\nu \Big)   \leq   \E(\phi( Y^1_T(\theta) + m) | \Fc_\nu)\\
 + \E \Big( \int_{\nu \vee \theta}^T \phi'(Y^1_s(\theta) + m) K_q (1 + \|Y^1 (\theta)\|_{\Sc^{\infty}[\theta,T]}) ds \Big| \Fc_\nu \Big) \;.
\end{multline*}
This last inequality and \reff{born-unif-Y} imply that there exists a constant $K$ which depends only on $m$, $T$ and $K_q$, such that for all $\F$-stopping times $\nu$ valued in $[0,T]$
\beqs
 \E \Big( \int_\nu^T| Z^1_s(\theta)|^2\mathds{1}_{s\geq\theta}  ds\Big|\Fc_\nu \Big) & \leq & K \;.
\enqs
\ep

\subsection{Some bounds about $X^0$ and $X^1$}
In this part, we give some bounds about the processes $X^0$ and $X^1$ which are used to get a uniform bound for the processes $Z^0$ and $Z^1$.
\begin{Proposition}\label{borne sur X} Suppose that \textbf{(HF)} holds. Then, we have
 \beq 
| \nabla X^0_t |~ :=~ \Big| \frac{\partial X^0_t}{\partial x} \Big| & \leq & e^{L_a T} \;, \quad 0\leq t\leq T\;,\label{majX0}
\enq
and  for any $\theta \in [0,T]$ we have
 \beq 
| \nabla^\theta X^1_t(\theta) | ~:=~ \Big| \frac{\partial X^1_t(\theta)}{\partial X^1_\theta(\theta)} \Big|  & \leq & e^{L_a T} \;, \quad \theta\leq t\leq T\;,\label{majX1}
\enq
 \beq 
| \nabla X^1_t(\theta) |~ :=~\Big| \frac{\partial X^1_t(\theta)}{\partial x} \Big|& \leq &(1 + L_a e^{L_a T}) e^{L_a T} \;,  \quad \theta \leq t\leq T\;.\label{majX1-x}
\enq
\end{Proposition}
\ni\textbf{Proof.} 
We first suppose that $b$ and $\beta$ are $C^1_b$ w.r.t. $x$. By definition we have 
 \beqs
\nabla X^0_t & = & 1 + \int_0^t \partial_x b ( s, X^0_s ) \nabla X^0_s ds\;,  \quad 0 \leq t \leq T \;.
\enqs 
  We get from Gronwall's lemma
 \beqs 
| \nabla X^0_t | & \leq & e^{L_a T} \;, \quad \quad0\leq t\leq T\;.
\enqs
In the same way, we have
\beqs 
\nabla^\theta X^1_t(\theta) & = & 1 + \int_\theta^t \partial_x b ( s, X^1_s(\theta)) \nabla^\theta X^1_s(\theta) ds \;, \quad \quad \theta \leq t \leq T \;,
\enqs 
and from Gronwall's lemma we get
 \beqs 
| \nabla^\theta X^1_t | & \leq & e^{L_a T} \;, \quad \theta\leq t\leq T\;.
\enqs
Finally we prove the last inequality. By definition 
\beqs 
\nabla X^1_t(\theta)  & = & 1 + \int_0^t \partial_x b ( s, X^1_s(\theta)) \nabla X^1_s(\theta) ds + \partial_x \beta (\theta, X^0_\theta) \nabla X^0_\theta \;, \quad \quad \theta \leq t \leq T \;.
\enqs 
Using the inequality \reff{majX0}, we get
\beqs 
|\nabla X^1_t(\theta)| & \leq & 1+ L_a e^{L_a T} + \int_0^t L_a |\nabla X^1_s(\theta)| ds  \;, \quad \quad \theta \leq t \leq T \;,
\enqs 
from Gronwall's lemma we get
 \beqs 
| \nabla X^1_t(\theta) | & \leq &(1 + L_a e^{L_a T}) e^{L_a T} \;,  \quad \theta \leq t\leq T\;.
\enqs
When $b$ and $\beta$ are not differentiable, we can also prove the result by regularization. We consider a density $q$ which is $C^\infty_b$ on $\R$ with a compact support, and we define an approximation $(b^\epsilon,\beta^\epsilon)$ of $(b, \beta)$ in $C^1_b$ by
\beqs
(b^\epsilon,  \beta^\epsilon)(t, x) & = & \frac{1}{\eps}\int_{\R} (b, \beta) (t, x') q \Big( \frac{x-x'}{\epsilon} \Big) dx'\;,\quad (t,x)\in[0,T]\times\R\;.
\enqs
We then use the convergence of $(X^{0,\epsilon},X^{1, \epsilon} (\theta))$ to $(X^0,X^{1} (\theta))$ and we get the result.
\ep

\subsection{Some bounds about $Y^0$ and $Y^1$}
In this part, we give some bounds about the processes $Y^0$ and $Y^1$ which are used to get a uniform bound for the processes $Z^0$ and $Z^1$.
\begin{Lemma}
Suppose that \textbf{(HF)}, \textbf{(HBQ)} and \textbf{(HBQD)} hold. Then, for any $\theta \in [0, T]$
\beq \label{majY1}
| \nabla^\theta Y^1_t(\theta)| ~:= ~\Big| \frac{\partial Y^1_t(\theta)}{\partial X^1_\theta(\theta)} \Big| & \leq & e^{( L_a + K_f ) T } ( K_g + T K_f ) \;,\quad \theta\leq t\leq T\;.
\enq 
\end{Lemma}

\ni\textbf{Proof.}
We first suppose that $b$, $f$ and $g$ are $C^1_b$ w.r.t. $x$, $y$ and $z$. In this case 
$(X^1(\theta) ,  Y^1(\theta), Z^1(\theta))$ is also differentiable w.r.t. $X^1_\theta(\theta) $ and we have
\beq \label{gradient de Y}
 \nabla^\theta Y^1_t(\theta) 
 & = & \nabla g ( X^1_T(\theta)  ) \nabla^\theta X^1_T(\theta)  - \int_t^T \nabla^\theta Z^1_s(\theta) dW_s \\
& & \nonumber + \int_t^T  \nabla f (s, X^1_s(\theta) , Y^1_s(\theta), Z^1_s(\theta) , 0 )  \big( \nabla^\theta X^1_s(\theta) ,\nabla^\theta Y^1_s(\theta) , \nabla^\theta Z^1_s(\theta) \big)  ds \;,
\enq
for $t\in[\theta,T]$. 
Define the process $R(\theta)$ by
\beqs
R_t(\theta) & := & \exp\Big({\int_0^t \partial_y  f(s, X^1_s(\theta), Y^1_s(\theta), Z^1_s(\theta),0)\mathds{1}_{s\geq \theta}ds}\Big)\;,\quad 0\leq t\leq T\;.
\enqs  
\ni Applying It\^o's formula, we get
\beq\nonumber
R_t(\theta) \nabla^\theta Y^1_t(\theta) & = & R_T(\theta) \nabla g(X^1_T(\theta)) \nabla^\theta X^1_T(\theta)  \\
& &  + \int_t^T R_s(\theta) \partial_x f (s, X^1_s(\theta), Y^1_s(\theta), Z^1_s(\theta),0 ) \nabla^\theta X^1_s(\theta) d s \nonumber\\
& & - \int_t^T R_s(\theta) \nabla^\theta Z^1_s(\theta) d W^1_s(\theta) \;, \quad \theta \leq t \leq T \;,\label{itoY1}
\enq 
where  the process $W^1(\theta)$ is defined by
\beq\label{W1}
 W^1_t(\theta) := W_t- \int_0^t\partial_z f (s, X^1_s(\theta), Y^1_s(\theta), Z^1_s(\theta), 0)\mathds{1}_{s\geq\theta} ds
\enq
for $t\in[0,T]$.
 From \textbf{(HBQD)}, there exists a constant $K >0$ such that we have 
\beqs
\Big\| \int_0^.\partial_z f (s, X^1_s(\theta), Y^1_s(\theta), Z^1_s(\theta), 0)\mathds{1}_{s\geq\theta} dW_s \Big\|^2_{BMO_\F[0,T]} & \leq  &\\
 K  \Big( 1+\sup_{\vartheta\in\Tc_\F[0,T]} \E\Big[\int_\vartheta^T|Z_s^1(\theta)|^2\mathds{1}_{s \geq \theta} ds\Big|\Fc_\vartheta\Big]\Big) 
 & \leq & \\
 K \Big(1+\Big\|\int_0^.Z^1_s(\theta)\mathds{1}_{ s\geq \theta}dW_s\Big\|^2_{BMO_\F[0,T]}\Big)
  & < & \infty\;,
\enqs
where the last inequality comes from Lemma \ref{BMO}. 

Hence by Lemma \ref{ppte-gen-BMO} the process $\Ec(\int_0^.\partial_z f (s, X^1_s(\theta), Y^1_s(\theta), Z^1_s(\theta), 0) \mathds{1}_{s\geq\theta}dW_s )$ is a uniformly integrable martingale.
Therefore, under the probability measure $\Q^1(\theta)$ defined by 
\beqs
\frac{d \Q^1(\theta)}{d \P} \Big|_{\Fc_t} & := & \Ec\Big(\int_0^.\partial_z f (s, X^1_s(\theta), Y^1_s(\theta), Z^1_s(\theta), 0) \mathds{1}_{s\geq\theta}dW_s \Big)_t\;,\quad 0\leq t\leq T\;,
\enqs
 we can apply Girsanov's theorem and $W ^1(\theta)$ is a Brownian motion under the probability measure $\Q^1(\theta)$. We then get from  \reff{itoY1}
\beqs
R_t(\theta) \nabla^\theta Y^1_t & = & \E_{\Q^1(\theta)}\Big[R_T(\theta) \nabla g(X^1_T) \nabla^\theta X^1_T  
 + \int_t^T R_s(\theta) \partial_x f (s, X^1_s, Y^1_s, Z^1_s,0 ) \nabla^\theta X^1_s d s \Big|\Fc_t\Big]\;.
\enqs
This last equality, \textbf{(HBQD)} and \reff{majX1}  give
\beq \label{majY1}
| \nabla^\theta Y^1_t(\theta)| & \leq & e^{( L_a + K_f ) T } ( K_g + T K_f ) \;,\quad \theta\leq t\leq T\;.
\enq 
 When $b$, $f$ and $g$ are not $C^1_b$, we can also prove the result by regularization as for Proposition \ref{borne sur X}.
\ep

\begin{Lemma}Suppose that \textbf{(HF)},  \textbf{(HBQ)} and \textbf{(HBQD)} hold. Then,
\beq
| \nabla Y^1_t(t)|  &\leq&  (1+L_a e^{L_a T})e^{( L_a + K_f ) T } ( K_g + T K_f )\;,  \quad 0 \leq t \leq T \;. \label{maj nabla Y1}
\enq
\end{Lemma}

 \ni\textbf{Proof.}
 Firstly, we suppose that $b$, $\beta$, $g$ and $f$ are $C^1_b$ w.r.t. $x$, $y$ and $z$. Then, for any $t \in [0,T]$
 \beq \label{gradient de Y1}
 \nabla Y^1_t(t) ~:= ~ \frac{\partial Y^1_t(t)}{\partial x } & = & \nabla g ( X^1_T(t)) \nabla X^1_T(t)\\
\nonumber & & + \int_t^T  \nabla f (s, X^1_s(t), Y^1_s(t), Z^1_s(t), 0)(\nabla X^1_s(t),\nabla Y^1_s(t),\nabla Z^1_s(t))ds\\
 & &  \nonumber   - \int_t^T \nabla Z^1_s(t) dW_s\;.
\enq
\ni Applying It\^o's formula, we get
\beqs\nonumber
R_t(t) \nabla Y^1_t(t) & = & R_T(t) \nabla g(X^1_T(t)) \nabla X^1_T(t)  \\
& &  + \int_t^T R_s(t) \partial_x f (s, X^1_s(t), Y^1_s(t), Z^1_s(t),0 ) \nabla X^1_s(t) d s \nonumber\\
& & - \int_t^T R_s(t) \nabla Z^1_s(t) d W^1_s(t) \;, \quad 0 \leq t \leq T \;,\label{itoY1-x}
\enqs 
where  the process $W^1(.)$ is defined in \reff{W1}. We have proved previously that  $W ^1(t)$ is a Brownian motion under the probability measure $\Q^1(t)$. We then get 
\beqs
R_t(t) \nabla Y^1_t(t) & = & \E_{\Q^1(t)}\Big[R_T(t) \nabla g(X^1_T(t)) \nabla X^1_T(t)  
 + \int_t^T R_s(t) \partial_x f (s, X^1_s(t), Y^1_s(t), Z^1_s(t),0 ) \nabla X^1_s(t) d s \Big|\Fc_t\Big]\;.
\enqs
This last inequality, \textbf{(HBQD)} and \reff{majX1}  give
\beqs
| \nabla Y^1_t(t)|  &\leq&  (1+L_a e^{L_a T})e^{( L_a + K_f ) T } ( K_g + T K_f )\;,  \quad 0 \leq t \leq T \;.
\enqs
When $b$, $f$ and $g$ are not $C^1_b$, we can also prove the result by regularization as for Proposition \ref{borne sur X}.
 \ep
 
\begin{Lemma}\label{nabla Y 0}
Suppose that \textbf{(HF)},  \textbf{(HBQ)} and \textbf{(HBQD)}  hold. Then, 
\beqs
| \nabla Y^0_t | & \leq & e^{ (2 K_f + L_a) T } (K_g + K_f T) \big( 1 + T K_f e^{K_fT} (1 + L_a e^{L_a T})\big) \;, \quad 0 \leq t \leq T\;.
\enqs 
\end{Lemma}
 
 \ni\textbf{Proof.}
We first suppose that $b$, $\beta$, $g$ and $f$ are $C^1_b$ w.r.t. $x$, $y$, $z$ and $u$, then $(X^0, Y^0, Z^0)$ is differentiable w.r.t. $x$ and  we have 
 \beqs \label{gradient de Y0}
\nabla Y^0_t& = & \nabla g ( X^0_T) \nabla X^0_T\\
\nonumber & & + \int_t^T \Big( \nabla f (s, X^0_s, Y^0_s, Z^0_s, Y^1_s(s) - Y^0_s)(\nabla X^0_s,\nabla Y^0_s,\nabla Z^0_s,\nabla Y^1_s(s)-\nabla Y^0_s)ds\\
  & &  \nonumber 
  - \int_t^T \nabla Z^0_s dW_s\;,  \quad 0 \leq t \leq T\;.
\enqs
Define the process $R^0$ by 
\beqs
R_t^0 & := & \exp\Big(\int_0^t (\partial_y - \partial_u) f (s, X^0_s, Y^0_s, Z^0_s, Y^1_s(s) - Y^0_s)ds\Big)\;,\quad 0\leq t\leq T\;.
\enqs
Applying It\^o's fomula we have
\beqs
R^0_t\nabla Y^0_t & = & R_T^0 \nabla g(X^0_T) \nabla X^0_T  \\
& &  + \int_t^T R_s^0 \partial_x f (s, X^0_s, Y^0_s, Z^0_s, Y^1_s(s) - Y^0_s ) \nabla X^0_s d s \\
& & + \int_t^T R_s^0   \partial_u f (s, X^0_s, Y^0_s, Z^0_s, Y^1_s(s) - Y^0_s ) \nabla Y^1_s ( s )  ds\\
& & - \int_t^T R_s^0  \nabla Z^0_s d {W}^0_s 
\enqs
where $d {W}^0_s := d W_s - \partial_z f (s, X^0_s, Y^0_s, Z^0_s, Y^1_s(s) - Y^0_s) ds$. From \textbf{(HBQD)}, there exists a constant $K >0$ such that we have 
\beqs
\Big\| \int_0^.\partial_z f (s, X^0_s, Y^0_s, Z^0_s, Y^1_s(s) - Y^0_s) dW_s \Big\|^2_{BMO_\F[0,T]} & \leq  & \\K  \Big( 1+\sup_{\vartheta\in\Tc_\F[0,T]} \E\Big[\int_\vartheta^T|Z_s^0|^2 ds\Big|\Fc_\vartheta\Big]\Big)
  &\leq& \\  K \Big(1+\Big\|\int_0^.Z^0_sdW_s\Big\|^2_{BMO_\F[0,T]}\Big)
  & < & \infty\;,
\enqs
where the last inequality comes from Lemma \ref{BMO}. 

Hence by Lemma \ref{ppte-gen-BMO}  the process $\Ec(\int_0^.\partial_z f (s, X^0_s, Y^0_s, Z^0_s, Y^1_s(s) - Y^0_s) dW_s )$ is a uniformly integrable martingale. Therefore, under the probability measure $\Q^0$ defined by
\beqs
\frac{d \Q^0}{d \P} \Big|_{\Fc_t} &  := &  \Ec\Big(\int_0^.\partial_z f (s, X^0_s, Y^0_s, Z^0_s, Y^1_s(s) - Y^0_s) dW_s\Big)_t
\enqs
 we can apply Girsanov's theorem and $W ^0$ is a Brownian motion under the probability measure $\Q^0$. Then, we get 
\beqs
R_t^0 \nabla Y^0_t & = & \E_{\Q^0} \Big[ R_T^0 \nabla g(X^0_T) \nabla X^0_T  \\
& &  + \int_t^T R_s^0 \partial_x f (s, X^0_s, Y^0_s, Z^0_s, Y^1_s(s) - Y^0_s ) \nabla X^0_s d s \\
& & + \int_t^T R_s^0 \partial_u f (s, X^0_s, Y^0_s, Z^0_s , Y^1_s(s) - Y^0_s) \nabla Y^1_s ( s )  ds \Big| \Fc_t \Big] \;.
\enqs
Using inequalities \reff{majX0} and \reff{maj nabla Y1} we get
\beqs
| \nabla Y^0_t | & \leq & e^{ (2 K_f + L_a) T } (K_g + K_f T) \big( 1 +  K_f T e^{K_fT} (1 + L_a e^{L_a T})\big) \;, \quad 0 \leq t \leq T\;.
\enqs
When $b$, $\beta$, $f$ and $g$ are not $C^1_b$, we can also prove the result by regularization as for Proposition \ref{borne sur X}.
\ep

\subsection{A uniform bound for $Z^0$ and $Z^1$}
Using the previous bounds, we obtain a uniform bound for the processes $Z^0$ and $Z^1$.
\begin{Proposition} \label{borne 1 pour Z 1}
Suppose that \textbf{(HF)}, \textbf{(HBQ)} and \textbf{(HBQD)} 
  hold. Then, for any $\theta \in [0, T]$, there exists a version of $Z^1 ( \theta )$ such that
\beqs
| Z^1_t ( \theta ) | & \leq &e^{(  2L_a + K_f ) T } ( K_g + T K_f ) K_a \;, \quad \theta \leq t \leq T \;.
\enqs 
\end{Proposition}

\ni\textbf{Proof.} 
Using Malliavin calculus, we have the classical representation of the process $Z^1(\theta)$ given by $\nabla^\theta Y^1(\theta)(\nabla^\theta X^1(\theta))^{-1}\sigma(.)$ (see Section \ref{lipschitz}). In the case where $b$, $f$ and $g$ are $C^1_b$ w.r.t. $x$, $y$ and $z$, we obtain from \reff{majY1}
\beqs 
| Z^1_t(\theta)  | & \leq &  e^{(  2L_a + K_f ) T } ( K_g + T K_f ) K_a \quad a.s.
\enqs
since $|(\nabla^\theta X^1(\theta))^{-1}| \leq e^{L_a T}$ (the proof of this inequality is similar to the one of  \reff{majX1}).\\
When $b$, $f$ and $g$ are not differentiable, we can also prove the result by a standard approximation and stability results for BSDEs with linear growth. 
\ep

\vspace{3mm}

\begin{Proposition}\label{borne 1 pour Z 0}
Suppose that \textbf{(HF)},  \textbf{(HBQ)} and \textbf{(HBQD)} hold. Then, there exists a version of $Z^0$ such that
\beqs
| Z^0_t | &\leq &e^{ 2 ( K_f + L_a) T } (K_g + K_f T) \big( 1 + T K_f e^{K_fT} (1 + L_a e^{L_a T})\big) K_a \;,  \quad 0 \leq t \leq T \;.
\enqs 
\end{Proposition}
 
 \ni\textbf{Proof.}
Thanks to the Malliavin calculus, it is classical to show that a version of $Z^0$ is given by $\nabla Y^0 (\nabla X^0)^{-1} \sigma(.)$ (see Section \ref{lipschitz}). So, in the case where $b$, $\beta$, $g$ and $f$ are $C^1_b$ w.r.t. $x$, $y$, $z$ and $u$, we obtain from \reff{majX0} and Lemma \ref{nabla Y 0}
\beqs
| Z^0_t | & \leq & e^{ 2 ( K_f + L_a) T } (K_g + K_f T) \big( 1 + T K_f e^{K_fT} (1 + L_a e^{L_a T})\big) M_\sigma \quad a.s. 
\enqs
since $| (\nabla X^0_t)^{-1} | \leq e^{L_a T}$ (the proof of this inequality is similar to the one of \reff{majX0}).

When $b$, $\beta$, $g$ and $f$ are not differentiable, we can also prove the result by a standard approximation and stability results for BSDEs with linear growth. 
 \ep

\subsection{Convergence of the scheme for the BSDE}

\begin{Theorem}
Under  \textbf{(HF)}, \textbf{(HFD)}, \textbf{(HBQ)} and \textbf{(HBQD)}  we have the following estimate
\beqs
\sup_{t\in[0,T] } \E \Big[ \big| Y_{t} - Y^{\pi}_{t} \big|^2 \Big]  
+ \E \Big[ \int_0^T \big| Z_t - Z^\pi_{t} \big|^2 dt \Big] + \E \Big[ \int_0^T \lambda_t \big| U_t - U^\pi_{t} \big|^2   dt \Big] & \leq & K |\pi| \;,
\enqs
for a constant $K$ which does not depend on $\pi$.
\end{Theorem}
\ni\textbf{Proof.} Fix $M\in\R$ such that  
\beqs
M & \geq &  \max\Big\{e^{(  2L_a + K_f ) T } ( K_g + T K_f ) K_a ~; \\
 & &   \quad \quad e^{ 2 ( K_f + L_a) T } (K_g + K_f T) \big( 1 + T K_f e^{K_fT} (1 + L_a e^{L_a T})\big) K_a \Big\}\;,
\enqs 
and define the function $\tilde f$ by 
\beqs
\tilde f (t,x,y,z,u) & = & f (t,x,y,\varphi_M(z),u) \;,\quad (t,x,y,z,u)\in[0,T]\times\R\times\R\times\R\times\R\;,
\enqs 
 where 
 \begin{equation*}
 \varphi_M(z)    := 
 \left\{
 \begin{array}{ccc}
 z & \mbox{ if }& |z|\leq M \\
M\frac{z}{|z|}    & \mbox{ if }& |z|> M
 \end{array}
 \right. , \quad z\in\R\;.
 \end{equation*}
 We notice that $\varphi_M$ is Lipschitz continuous and bounded. Therefore  we obtain from  $\textbf{(HBQD)}$ that $\tilde f$ is Lipschitz continuous.  
 
 Moreover, using Propositions \ref{borne 1 pour Z 1} and \ref{borne 1 pour Z 0}, we get that under 
 \textbf{(HF)}, \textbf{(HBQ)} and \textbf{(HBQD)}, 
 $(X,Y,Z)$ is also solution to the Lipschitz FBSDE
 \beqs
  X_t  & = & x + \int_0^t b(s, X_s) ds + \int_0^t \sigma(s,X_s) dW_s + \int_0^t \beta(s, X_{s^-} ) dH_s \;, \quad 0 \leq t\leq T\;,\\ \nonumber
  Y_t  & =  & g(X_T) + \int_t^T\tilde  f\big(s, X_s, Y_s, Z_s, U_s(1-H_s)\big) ds\\ 
   & & \quad \quad \quad \quad  - \int_t^T Z_s dW_s - \int_t^T U_s dH_s\;, \quad 0 \leq t\leq T\;. ~~
 \enqs
 Applying Theorem \ref{erreur schema}, we get the result.
\ep

\end{document}